\documentclass[preprint,3p,12pt]{elsarticle}
\usepackage{mathrsfs}
\usepackage{amsmath}
\usepackage{stmaryrd}
\usepackage{bbding}
\usepackage{dcolumn}
\usepackage{graphicx}
\usepackage{amsfonts}
\usepackage{amssymb}
\usepackage{psfrag}
\usepackage{wrapfig}
\usepackage{subfigure}
\usepackage{makeidx}
\usepackage{bm}
\usepackage{epsf}
\usepackage{epsfig}
\usepackage{setspace}
\usepackage{graphicx}
\usepackage{epstopdf}
\usepackage{psfrag}
\usepackage{subfigure}
\usepackage{color}

\begin{document}

\title{Performance evaluation of high-order compact and second-order gas-kinetic schemes in compressible flow simulations }

\author[HKUST1]{Yaqing Yang}
\ead{yangyq@ust.hk}  

\author[HKUST1]{Fengxiang Zhao}
\ead{fzhaoac@connect.ust.hk}

\author[HKUST1,HKUST2]{Kun Xu \corref{cor}}
\ead{makxu@ust.hk}

\cortext[cor]{Corresponding author}
\address[HKUST1]{Department of Mathematics, Hong Kong University of Science and Technology, Clear Water Bay, Kowloon, Hong Kong}
\address[HKUST2]{Shenzhen Research Institute, Hong Kong University of Science and Technology, Shenzhen, China}
 
\begin{abstract}

The trade-off among accuracy, robustness, and computational cost remains a key challenge in simulating complex flows.
Second-order schemes are computationally efficient but lack the accuracy required for resolving intricate flow structures, particularly in turbulence. High-order schemes, especially compact high-order schemes, offer superior accuracy and resolution at a relatively modest computational cost.
To clarify the practical performance of high-order schemes in scale-resolving simulations, this study evaluates two representative gas-kinetic schemes: the newly developed fifth-order compact gas-kinetic scheme (CGKS-5th) and the conventional second-order gas-kinetic scheme (GKS-2nd).
Test cases ranging from subsonic to supersonic flows are used to quantitatively assess their accuracy and efficiency.
The results demonstrate that CGKS-5th achieves comparable resolution to GKS-2nd at roughly an order of magnitude lower computational cost. 
Under equivalent computational resources, CGKS-5th delivers significantly higher accuracy and resolution, particularly in turbulent flows involving shocks and small-scale vortices. 
This study provides the first clear verification of the advantages of high-order compact gas-kinetic schemes in simulating viscous flows with discontinuities.
Additionally, multi-GPU parallelization using CUDA and MPI is implemented to enable large-scale applications.

\end{abstract}
 

\maketitle

\section{Introduction}

To address complex flow simulations in engineering applications, it is essential to balance resolution, robustness, and computational cost. 
Most commercial and industrial software currently relies on second-order numerical methods, such as MUSCL \cite{MUSCL}, owing to their robustness, low computational cost, and broad applicability. However, second-order schemes are inherently constrained by significant numerical dissipation and dispersion, which limit their accuracy in high-resolution scenarios. For example, accurately predicting the position and intensity of shock waves in high-speed flows, capturing wave propagation and reflection in acoustic problems, and resolving multi-scale turbulence structures in direct numerical simulations (DNS) all expose the limitations of second-order methods.

To overcome these challenges, high-order schemes have been developed to provide improved accuracy and resolution for complex flow problems. Over the past few decades, various high-order methods have been proposed, including finite difference (FD) methods \cite{FD-1,FD-2}, essential non-oscillatory (ENO) schemes \cite{ENO-1,ENO-2}, weighted essential non-oscillatory (WENO) schemes \cite{WENO-1,WENO-2}, Hermite WENO (HWENO) schemes \cite{HWENO-1}, discontinuous Galerkin (DG) methods \cite{DG-0,DG-1}, flux reconstruction (FR) methods \cite{FR}, correction procedure using reconstruction (CPR) methods \cite{CPR}, and compact gas-kinetic scheme (CGKS) \cite{CGKS-high-2,CGKS-high-3,CGKS-high-4,CGKS-high-6,CGKS5-prep}. Finite difference and finite volume methods typically achieve high-order accuracy by using expanded stencils that incorporate additional neighboring information for reconstruction. In contrast, DG-type methods improve accuracy by increasing the degree of interpolation polynomials within each cell, thereby enhancing the internal degrees of freedom and naturally exhibiting compactness. However, the widespread adoption of high-order schemes in industrial engineering is limited by lower robustness, higher computational costs, and greater implementation complexity. Achieving an optimal balance among resolution, robustness, and computational efficiency is essential to fully exploit the potential of high-order schemes in practical applications.

In recent decades, the gas-kinetic scheme (GKS) has been systematically developed for simulating flows across a wide range of regimes, from subsonic to supersonic \cite{GKS-Xu1,GKS-Xu2,GKS-Xu3}. GKS performs a one-step, second-order-accurate temporal evolution from kinetic to hydrodynamic scales and offers several intrinsic advantages: it computes both inviscid and viscous fluxes within a single unified formulation and adaptively transitions between equilibrium-state fluxes in smooth regions and non-equilibrium transport fluxes suitable for discontinuities. In conjunction with linear reconstruction and slope limiters, the conventional second-order GKS can be obtained.
Building on this foundation, the high-order compact GKS (CGKS) further exploits the gas distribution function to achieve time-accurate evolution of flow variables at cell interfaces. As a result, both the cell-averaged conservative variables and their gradients within each control volume can be updated simultaneously at the next time level. Coupled with the two-stage temporal discretization method \cite{s2o4-0,GRP-high-1,s2o4-1}, high-order CGKS have been developed for high-fidelity simulation of compressible flows \cite{CGKS-high-2,CGKS-high-3,CGKS-high-4,CGKS-high-6,CGKS5-prep}.
For structured meshes, a high-order CGKS was proposed in \cite{CGKS-high-2}, demonstrating spectral-like resolution in one dimension using a dimension-by-dimension reconstruction strategy. However, this approach substantially increases algorithmic complexity in three-dimensional applications. High-order GKS in curvilinear coordinates have been developed using dimension-by-dimension reconstruction for both laminar \cite{HGKS-curv-1} and turbulent flows \cite{HGKS-curv-2}. Nevertheless, when extended to compact schemes in three dimensions, the dimension-by-dimension strategy still introduces significant implementation challenges.
To achieve high resolution and efficiency while maintaining implementation simplicity, a multidimensional fifth-order CGKS has been newly developed for non-orthogonal structured meshes \cite{CGKS5-prep}. This method employs line-averaged derivatives within each control volume to provide additional degrees of freedom. Its resolution, robustness, and stability have been validated through a series of complex flow simulations.

In this paper, we conduct a systematic performance evaluation of the newly developed fifth-order compact gas-kinetic scheme (CGKS-5th) against the conventional second-order gas-kinetic scheme (GKS-2nd) on structured meshes. Two comparative strategies are employed: (i) assessing the computational cost required to achieve a target resolution, and (ii) evaluating accuracy and resolution under comparable computational resources. 
While the advantages of high-order schemes over second-order schemes have been extensively validated in complex smooth flows \cite{High-perf-1,High-perf-2,High-perf-3}, such verification remains limited for flows involving discontinuities. Theoretically, the use of nonlinear reconstructions or limiters in discontinuous flows may impact the performance of high-order schemes. This study seeks to address this gap.
To achieve this, a series of numerical tests, ranging from subsonic to supersonic turbulent flows, are conducted to compare the resolution and efficiency of the two schemes under these criteria.
Furthermore, both schemes are implemented on multiple graphics processing units (GPUs) using the Compute Unified Device Architecture (CUDA) and the Message Passing Interface (MPI). Numerical experiments are conducted on NVIDIA GeForce RTX 4090 GPUs for achieving significant computational acceleration.

This paper is organized as follows: Section 2 introduces the gas-kinetic scheme and finite volume framework. Section 3 provides a brief review of the new fifth-order compact gas-kinetic scheme and the second-order gas-kinetic scheme. Numerical examples are presented in Section 4, and conclusions are drawn in the final section.

\section{Time-accurate solution of GKS and finite volume method}

In the gas-kinetic method, time-dependent numerical fluxes are computed from a physically modeled, space-time coupled solution of the BGK equation for the gas distribution function \cite{GKS-Xu1,GKS-Xu2}.
In the gas-kinetic method, a physically modeled, space-time coupled solution of the BGK equation for the gas distribution function is employed to compute time-dependent numerical fluxes \cite{GKS-Xu1,GKS-Xu2}.
The BGK equation is a simplification of Boltzmann equation, and the three-dimensional BGK equation \cite{BGK-1,BGK-2} can be written as
\begin{equation}\label{bgk}
f_t+\boldsymbol{u}\cdot \nabla f=\frac{g-f}{\tau},
\end{equation}
where $\boldsymbol{u}=(u,v,w)$ is the particle velocity, $\tau$ is the collision time and $f$ is the gas distribution function. The
equilibrium state $g$ is given by Maxwellian distribution 
\begin{equation*}
g=\rho\left(\frac{\lambda}{\pi}\right)^{(N+3)/2}e^{-\lambda[(\boldsymbol{u}-\boldsymbol{U})^2+\boldsymbol{\xi}^2]},
\end{equation*}
where $\rho$ is the density, $\boldsymbol{U}=(U,V,W)$ is the macroscopic fluid velocity, and $\lambda=\displaystyle\frac{\rho}{2p}$, $p$ is the pressure. 
In the BGK equation, the collision operator involves a simple relaxation from $f$ to the local equilibrium state $g$. 
The variable $\boldsymbol{\xi}$ accounts for the internal degree of freedom of molecular motion, $\boldsymbol{\xi}^2=\xi_1^2+\dots+\xi_N^2$, where $N=(5-3\gamma)/(\gamma-1)$ is the number of
internal degree of freedom and $\gamma$ is the specific heat ratio. 
The collision term satisfies the compatibility condition
\begin{equation*}
\int \frac{g-f}{\tau}\boldsymbol{\psi} \, \text{d}\Xi=0,
\end{equation*}
where $\displaystyle\boldsymbol{\psi}=\left(1,u,v,w,\frac{1}{2}(u^2+v^2+w^2+\boldsymbol{\xi}^2)\right)^T$ and $\text{d}\Xi=\text{d}u\text{d}v\text{d}w\text{d}\xi_1\dots\text{d}\xi_{N}$.
The macroscopic conservative variables $\boldsymbol{Q}=(\rho, \rho U,\rho V, \rho W, \rho E)$ can be calculated through the gas distribution function $f$
\begin{equation}\label{Q-macro}
\boldsymbol{Q}=\int f \boldsymbol{\psi} \,\text{d}\Xi,
\end{equation}
and the corresponding fluxes can be given by taking moments of the gas
distribution function
\begin{align}\label{flux-macro}
\boldsymbol{F}=\int \boldsymbol{u} f \boldsymbol{\psi} \, \text{d}\Xi,
\end{align}
where $\boldsymbol{F}$ represents either the Euler flux or the Navier-Stokes (NS) flux, depending on the order of approximation of $f$ to $g$. This study focuses solely on viscous flows described by the NS equations \cite{GKS-Xu1,GKS-Xu2}.
In GKS, $\boldsymbol{F}$ at cell interface is determined by the gas distribution function $f$.
Based on the integral solution of BGK equation 
Eq.\eqref{bgk}, the gas distribution function $f(\boldsymbol{x}_{G},t,\boldsymbol{u},\boldsymbol{\xi})$ can be given by
\begin{equation*}
f(\boldsymbol{x}_{G},t,\boldsymbol{u},\boldsymbol{\xi})=\frac{1}{\tau}\int_0^t
g(\boldsymbol{x}',t',\boldsymbol{u}, \boldsymbol{\xi})e^{-(t-t')/\tau}\text{d}t'+e^{-t/\tau}f_0(-\boldsymbol{u}t,\boldsymbol{\xi}).
\end{equation*}
By modeling and approximating the unknown terms in the integral solution, a second-order accurate explicit solution for the distribution function is obtained as
\begin{align}\label{flux}
f(\boldsymbol{x}_{G},t,\boldsymbol{u},\boldsymbol{\xi})=&(1-e^{-t/\tau})g_0+
[(t+\tau)e^{-t/\tau}-\tau](\overline{a}_1u+\overline{a}_2v+\overline{a}_3w)g_0\nonumber\\
+&(t-\tau+\tau e^{-t/\tau}){\bar{A}} g_0\nonumber\\
+&e^{-t/\tau}g_r[1-(\tau+t)(a_{1}^{r}u+a_{2}^{r}v+a_{3}^{r}w)-\tau A^r](1-H(u))\nonumber\\
+&e^{-t/\tau}g_l[1-(\tau+t)(a_{1}^{l}u+a_{2}^{l}v+a_{3}^{l}w)-\tau A^l]H(u),
\end{align}
where $g_{0}$, $g_l$ and $g_r$ are determined by the conservative variables reconstructed at the interface, and the coefficients are determined by the macroscopic variables and their derivatives. The specific formulas can be referred to in \cite{GKS-Xu1,GKS-Xu2}.

For the second-order evolution solution, the time accurate distribution function $f(t)$ at cell interface can be approximated through a linearization in time \cite{CGKS-high-5}
\begin{equation*}
\hat{f}(t)=f^n+t f_t^n.
\end{equation*}
The two coefficients $f^n$ and $f_t^n$ are calculated as follows 
\begin{align*}
    f^n&=\big(4\bar{f}(\Delta t/2) - \bar{f}(\Delta t)\big)/\Delta t,\\
    f_t^n&=4\big(\bar{f}(\Delta t) - 2\bar{f}(\Delta t/2)\big)/{\Delta t}^2,
\end{align*}
where $\bar{f}(\Delta t)$ and $\bar{f}(\Delta t/2)$ are the time integrations of $f(t)$ over the interval $[t_n, t_n + \Delta t]$ and $[t_n, t_n + \Delta t/2]$, respectively.
The numerical fluxes and their time derivatives can be obtained by taking moments of $\hat{f}(t)$ and $\hat{f}_t(t)$ at $t = t_n$
\begin{equation}\label{deri-t-f}
\boldsymbol{F}^n = \int \boldsymbol{u} f^n \boldsymbol{\psi} \, \text{d}\Xi,\quad
\boldsymbol{F}_t^n = \int \boldsymbol{u} f_t^n \boldsymbol{\psi} \, \text{d}\Xi.
\end{equation}
Simultaneously, the flow variables and their time derivatives can be given by
\begin{equation}\label{deri-t-q}
\boldsymbol{Q}^n = \int f^n \boldsymbol{\psi} \, \text{d}\Xi,\quad
\boldsymbol{Q}_t^n = \int f_t^n \boldsymbol{\psi} \, \text{d}\Xi.
\end{equation}
The obtained conservative variables at interfaces are used to achieve compact reconstructions in high-order compact GKS.

Taking moments of Eq.\eqref{bgk} and integrating with respect to space, the semi-discretized finite volume scheme can be obtained as
\begin{align}\label{semi}
\frac{\text{d} \boldsymbol{Q}_i}{\text{d} t}=\mathcal{L}_{i}(\boldsymbol{Q}),
\end{align}
where $\boldsymbol{Q}_i$ is the cell-averaged conservative variables over cell $\Omega_{i}$. 
The operator $\mathcal{L}_{i}$ is defined as
\begin{equation*}
\mathcal{L}_{i}(\boldsymbol{Q})
=-\frac{1}{|\Omega_{i}|} \int_{\partial\Omega_{i}}\boldsymbol{F}(t)\boldsymbol{n}\,\text{d}S,
\end{equation*}
where $|\Omega_{i}|$ is the volume of $\Omega_{i}$, $\partial\Omega_{i}$ represents the cell interfaces of $\Omega_{i}$, and $\boldsymbol{n}$ is the unit direction of the cell interface. 
To achieve the expected order of accuracy, the Gaussian 
quadrature is used for the flux integration
\begin{align*}
\int_{\sigma_{j}}\boldsymbol{F}(t)\boldsymbol{n}\,\text{d}S=\sum_{k=1}^{K}\omega_{k}\boldsymbol{F}(\boldsymbol{x}_k,t)S_{j},
\end{align*}
where $S_{j}$ is the area of $\sigma_{j}$, $\sigma_{j}$ is one of the cell interfaces, $\omega_{k}$ is the Gaussian quadrature weight, $\boldsymbol{x}_{k}$ is the Gaussian quadrature point, $K$ is the total number of Gaussian quadrature points.
The numerical flux $\boldsymbol{F}(\boldsymbol{x}_k,t)$ at Gaussian quadrature point can be given by Eq.\eqref{flux-macro}.
In this study, structured meshes are utilized, where each face is a quadrilateral with $K=4$ Gaussian points for high-order schemes.
For the second-order scheme, this integral reduces to a form that only involves the midpoints of the cell interfaces.

\section{High-order compact and second-order gas-kinetic schemes}

Second-order schemes have been widely used to solve a broad range of engineering problems owing to their robustness and computational efficiency. However, their accuracy is often inadequate for accurately capturing complex flow features in high-fidelity simulations. In contrast, high-order schemes, due to their superior accuracy, can faithfully resolve intricate details of the flow field, although this comes at the cost of increased algorithmic complexity.
To enable a clear and fair assessment of the differences between high- and low-order methods in simulations, this study conducts a performance evaluation within the framework of the finite-volume method based on GKS. Specifically, we consider a recently proposed CGKS-5th \cite{CGKS5-prep} and a conventional GKS-2nd, as representatives of high-order and second-order approaches, respectively. A brief introduction to these two schemes is provided in this section.

\subsection{Flow variables update}

For both CGKS-5th and GKS-2nd, the cell-averaged conservative variable $\boldsymbol{Q}(t)$ within a control volume $\Omega$ is updated according to the conservation laws, as given in Eq.\eqref{flux}. In CGKS-5th, additionally, the cell-averaged and line-averaged derivatives of $\boldsymbol{Q}(t)$ are also updated.

Using Gauss's theorem, the cell-averaged gradient of $\boldsymbol{Q}(t)$ over $\Omega$ can be expressed as
\begin{equation*}\label{derivative}
\displaystyle\nabla \boldsymbol{Q}(t)=\frac{1}{|\Omega|}\int_{\Omega}\nabla \boldsymbol{Q}(t)\text{d}V=\frac{1}{|\Omega|}\int_{\partial\Omega}\boldsymbol{Q}(t) \boldsymbol{n}\text{d}S,
\end{equation*}
where $\boldsymbol{n}$ is the unit normal vector on the cell interface.
The line-averaged partial derivative along a segment defined by points $\boldsymbol{x}_{1}$ and $\boldsymbol{x}_{2}$ is given by
\begin{align*}
(\partial_{l} \boldsymbol{Q})(t) &= \frac{1}{\left|l_G\right|}\int_{\boldsymbol{x}_{1}}^{\boldsymbol{x}_{2}}\frac{\partial \boldsymbol{Q}(t)}{\partial l}\text{d}l \\
&= \frac{1}{\left|l_G\right|}\left(\boldsymbol{Q}(\boldsymbol{x}_{2},t)-\boldsymbol{Q}(\boldsymbol{x}_{1},t)\right),
\end{align*}
where $\left|l_G\right| = \left|\boldsymbol{x}{2}-\boldsymbol{x}{1}\right|$, and $\boldsymbol{Q}(\boldsymbol{x}_{1},t)$ and $\boldsymbol{Q}(\boldsymbol{x}_{2},t)$ are the values of $\boldsymbol{Q}$ at $\boldsymbol{x}_{1}$ and $\boldsymbol{x}_{2}$ inside the cell, respectively.
At each time step, $\boldsymbol{Q}(\boldsymbol{x}_{1},t)$ and $\boldsymbol{Q}(\boldsymbol{x}_{2},t)$ are obtained from the time-accurate solution of the gas distribution function in the GKS. This solution is also used to compute the numerical flux at the interfaces, and it is specified using interface values provided by spatial reconstruction.

\begin{figure}[!h]
\centering
\includegraphics[width=0.4\textwidth]{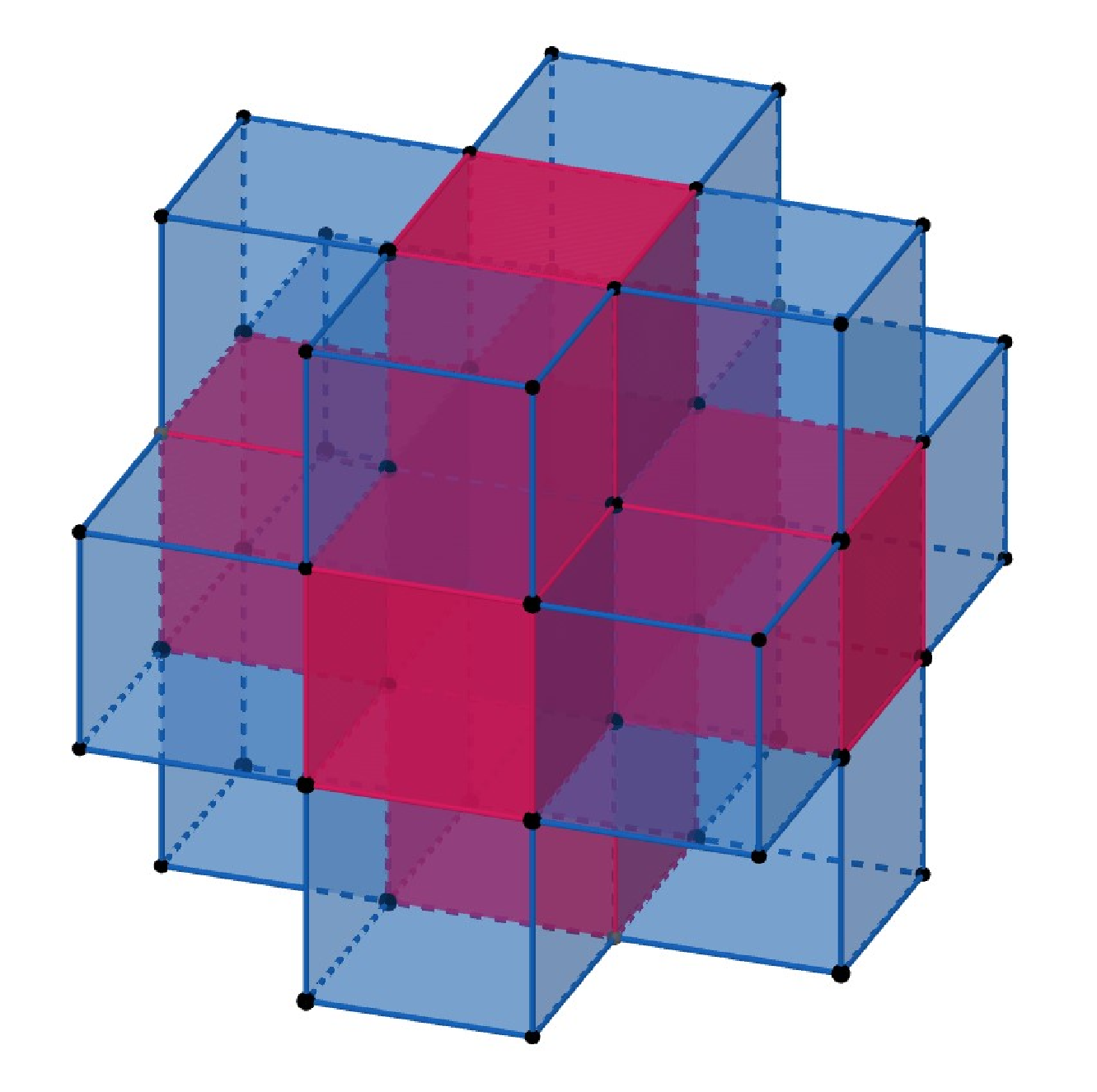}{(a)}
\includegraphics[width=0.4\textwidth]{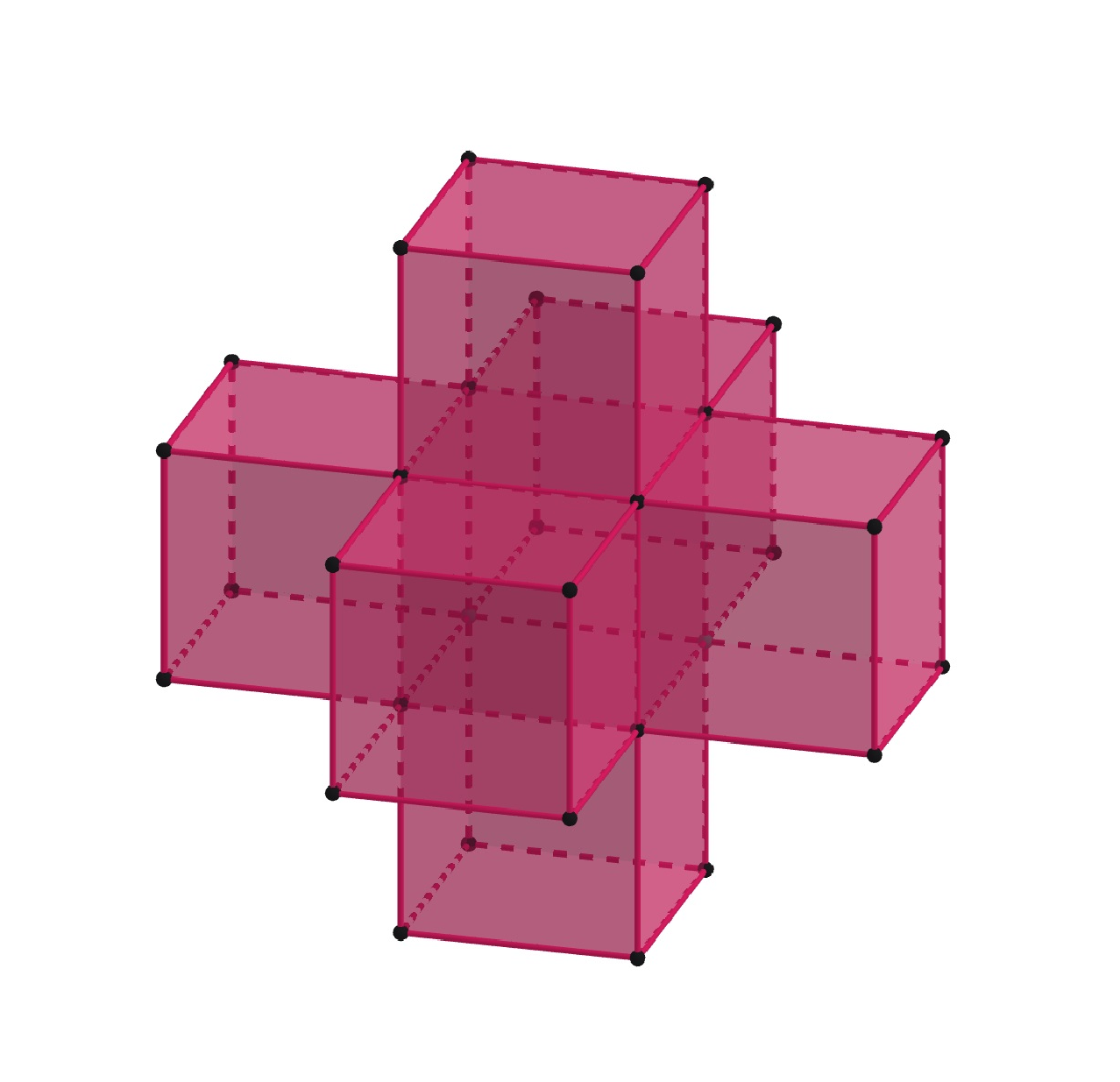}{(b)}
\caption{\label{stencil} The schematic diagram of the compact stencil used for spatial reconstruction in (a) CGKS-5th and (b) GKS-2nd.}
\end{figure}

\subsection{Fifth-order compact gas-kinetic scheme}\label{CGKS-5}

Figure \ref{stencil}(a) and Figure \ref{stencil}(b) illustrate the schematic of the compact stencils used for spatial reconstruction in CGKS-5th and GKS-2nd, respectively. In these figures, the red cubes denote face-neighboring cells, while the blue cubes indicate edge-neighboring cells.
For the target cell $\Omega_0$, its interfaces are denoted by $F_m$ ($m = 1, \dots, 6$), and its edges are denoted by $E_n$ ($n = 1, \dots, 12$). 
The cell adjacent to $\Omega_0$ that shares the face $F_m$ is represented by $\Omega_{1,m}$, while the cell adjacent to $\Omega_0$ that shares only the edge $E_n$ is represented by $\Omega_{2,n}$.

For CGKS-5th, to achieve fifth-order accuracy, a compact stencil containing all the face-neighboring and edge-neighboring cells is defined as follows
\begin{equation*}
S^{cell-5th}=\{\Omega_0, \Omega_{1,m}, \Omega_{2,n}\,|\,m = 1, \dots, 6,\, n = 1, \dots, 12\}.
\end{equation*}
The following data are utilized for fifth-order reconstruction
\begin{align*}\label{stencil-5}
S^{data-5th}=S^{0}\cup S^{face}\cup S^{edge},
\end{align*}
where
\begin{align*}
    S^{0}&=\{\boldsymbol{Q}_0,\,(\partial_l \boldsymbol{Q})_{0,G}\},\\
    S^{face}&=\{\boldsymbol{Q}_{1,m},\,(\nabla \boldsymbol{Q})_{1,m}\,|\,m = 1, \dots, 6\},\\
    S^{edge}&=\{\boldsymbol{Q}_{2,n},\,(\partial_{\boldsymbol{\tau}} \boldsymbol{Q})_{2,n}\,|\,n=1,\dots,12\}.
\end{align*}
$\boldsymbol{Q}_0$, $\boldsymbol{Q}_{1,m}$ and $\boldsymbol{Q}_{2,n}$ represent the cell-averaged conservative values over cells $\Omega_0$, $\Omega_{1,m}$ and $\Omega_{2,n}$, respectively. $(\partial_{l} \boldsymbol{Q})_{0,G}$ represents the line-averaged partial derivatives within the target cell $\Omega_0$. 
$(\nabla \boldsymbol{Q})_{1,m}$ denotes the cell-averaged gradient over cell $\Omega_{1,m}$, where $\nabla=(\partial_x,\partial_y,\partial_z)$. 
$(\partial_{\boldsymbol{\tau}} \boldsymbol{Q})_{2,n}$ denotes the cell-averaged directional derivative over the cell $\Omega_{2,n}$, with $\partial_{\boldsymbol{\tau}} = (\partial_{\tau_1}, \partial_{\tau_2})$, where $\partial_{\tau_1}$ and $\partial_{\tau_2}$ are linear combinations of $\partial_x$, $\partial_y$, and $\partial_z$. Further details on the computation of these quantities can be found in \cite{CGKS5-prep}.

For the target cell $\Omega_0$, a quartic polynomial $P^4(\boldsymbol{x})$ can be reconstructed based on $S^{data-5th}$, defined by
\begin{equation*}
P^4(\boldsymbol{x})=\boldsymbol{Q}_{0}+\sum_{|\boldsymbol{d}|=1}^4 a_{\boldsymbol{d}}p_{\boldsymbol{d}}(\boldsymbol{x}),
\end{equation*}
where $\boldsymbol{Q}_{0}$ is the cell-averaged variable over cell $\Omega_{0}$, the multi-index $\boldsymbol{d}=(d_1, d_2, d_3)$, and $|\boldsymbol{d}|=d_1+d_2+d_3$. 
The zero-mean basis $p_{\boldsymbol{d}}(\boldsymbol{x})$ is defined as
\begin{equation}\label{base}
p_{\boldsymbol d}(\boldsymbol{x})=\displaystyle\frac{1}{d_1!d_2!d_3!}{\delta x}^{d_1}{\delta y}^{d_2}{\delta z}^{d_3}-\frac{1}{\left|\Omega_{0}\right|}\int_{\Omega_{0}}\frac{1}{d_1!d_2!d_3!}{\delta x}^{d_1}{\delta y}^{d_2}{\delta z}^{d_3}\text{d}V,
\end{equation}
where
\begin{equation*}
\delta x=\frac{1}{h_x}(x-x_0),\,\delta y=\frac{1}{h_y}(y-y_0),\delta z=\frac{1}{h_z}(z-z_0),
\end{equation*}
$h_x$, $h_y$ and $h_z$ are the characteristic scales of $\Omega_0$ along the three directions of axes, and $(x_0,y_0,z_0)$ is the centroid of $\Omega_0$.
To determine the polynomial $P^4(\boldsymbol{x})$, the following system is solved using the constrained least-squares method
\begin{equation}\label{compact-big-5}
\begin{split}
\frac{1}{\left|\Omega_{k}\right|}\int_{\Omega_{k}}P^4(\boldsymbol{x})\text{d}V&=\boldsymbol{Q}_{k},~\Omega_{k}\in S^{cell-5th},\\
\frac{h_1}{\left|\Omega_{1,m}\right|}\int_{\Omega_{1,m}}\nabla P^4(\boldsymbol{x})\text{d}V&=h_1\cdot(\nabla \boldsymbol{Q})_{1,m},~m=1,\dots, 6,\\
\frac{h_2}{\left|\Omega_{2,n}\right|}\int_{\Omega_{2,n}}\frac{\partial}{\partial_{\boldsymbol{\tau}}}P^4(\boldsymbol{x})\text{d}V&=h_2\cdot(\partial_{\boldsymbol{\tau}} \boldsymbol{Q})_{2,n},~n=1,\dots,12,\\
\frac{h_3}{\left|l_G\right|}\int_{l_G}\frac{\partial}{\partial_{l}}P^4(\boldsymbol{x})\text{d}l&=h_3\cdot(\partial_{l} \boldsymbol{Q})_{0,G},
\end{split}
\end{equation}
where the parameters $h_1$, $h_2$, and $h_3$ are used to reduce the condition number of the reconstruction matrix, the choice of which can be found in \cite{CGKS5-prep}. In the constrained least-squares system, the first-row equation in Eq.\eqref{compact-big-5} is satisfied exactly, while the remaining equations are satisfied in the least-squares sense.

To handle discontinuities, nonlinear reconstruction is necessary. For CGKS-5th, we adopt the generalized ENO (GENO) nonlinear reconstruction method \cite{GENO} to obtain nonlinear reconstruction. At discontinuities, the GENO reconstruction is dominated by a second-order reconstruction with ENO property. The second-order reconstruction is determined based on sub-stencils.
The sub-stencils $S_{m}^{cell-2}$ ($m = 1, \dots, 6$) of the CGKS-5th scheme are selected as
\begin{equation*}
S_{m}^{cell-2}=\{\Omega_0, \Omega_{1,m}\}.
\end{equation*}
The construction of several linear polynomials on these sub-stencils is based on the following two types of data
\begin{align*}\label{stencil-2}
S_{m}^{data-2}=\{\boldsymbol{Q}_0,\,(\partial_l \boldsymbol{Q})^m_{0,G}, \boldsymbol{Q}_{1,m}\},\,m=1,\dots,6,
\end{align*}
where
$\boldsymbol{Q}_{1,m}$ is the cell-averaged conservative value over cell $\Omega_{1,m}$, $(\partial_l \boldsymbol{Q})^m_{0,G}$ is the selected line-averaged partial derivative within the target cell $\Omega_0$. Detailed descriptions of this selection can be found in \cite{CGKS5-prep}.
Using $S_{m}^{data-2}$ ($m=1,\dots,6$), linear polynomials $P_m^1(\boldsymbol{x})$ can be constructed of the form
\begin{equation}\label{p1}
P_m(\boldsymbol{x})=\boldsymbol{Q}_{0}+\sum_{|\boldsymbol d|=1}b_{\boldsymbol d}^mp_{\boldsymbol d}(\boldsymbol{x}),
\end{equation}
where the basis functions $p_{\boldsymbol d}(\boldsymbol{x})$ are defined as in Eq.\eqref{base}.
To determine these linear polynomials, the following constrains need to be satisfied
\begin{equation}\label{compact-sub-2}
\begin{split}
\frac{1}{\left|\Omega_{1,m}\right|}\int_{\Omega_{1,m}}P_m^1(\boldsymbol{x})\text{d}V&=\boldsymbol{Q}_{1,m},\\
\frac{1}{\left|l_G\right|}\int_{l_G}\frac{\partial}{\partial_{l}}P_m^1(\boldsymbol{x})\text{d}l&=(\partial_{l} \boldsymbol{Q})^m_{0,G},
\end{split}
\end{equation}
where $\boldsymbol{Q}_{1,m}$ is the cell-averaged variable over cell $\Omega_{1,m}$.  The resulting linear system is solved using a constrained least-squares method.

To accommodate non-uniform meshes while keeping memory usage low and the implementation simple, the CGKS-5th scheme adopts a simple transformation strategy \cite{CGKS5-prep}. Specifically, a coordinate mapping sends each target cell to a single reference cell, allowing the reconstruction to be represented in a unified polynomial form. Consequently, only one reconstructed polynomial on the reference cell needs to be stored, reconstructed values are then mapped back from the reference coordinates to the physical coordinates.

Using reconstructions on the compact stencil and sub-stencils, GENO is applied to achieve nonlinear reconstruction. It adaptively transitions from high-order linear reconstruction in smooth regions to second-order reconstruction near discontinuities, effectively balancing accuracy and robustness.
With the reconstructed polynomial $P^4(\boldsymbol{x})$ and $P_m^1(\boldsymbol{x}), m=1,...,6$, the point-value $\boldsymbol{Q}(\boldsymbol{x}_{G})$ and the spatial derivatives $\partial_{x,y,z} \boldsymbol{Q}(\boldsymbol{x}_{G})$ at Gaussian quadrature point for the nonlinear reconstruction are computed as
\begin{equation}\label{weno-new}
\begin{split}
\boldsymbol{Q}(\boldsymbol{x}_{G})&=\chi P^4(\boldsymbol{x}_{G}) + (1-\chi)\left(\sum_{m=1}^{6}\omega_{m} P_m^1(\boldsymbol{x})\right),\\
\partial_{x,y,z} \boldsymbol{Q}(\boldsymbol{x}_{G})&=\chi \partial_{x,y,z} P^4(\boldsymbol{x}_{G}) + (1-\chi)\left(\sum_{m=1}^{6}\omega_{m} \partial_{x,y,z} P_m^1(\boldsymbol{x})\right),
\end{split}
\end{equation}
where 
\begin{equation*}
\displaystyle\omega_{m}=\frac{\overline{\omega_{m}}}{\sum_{m=1}^{6} \overline{\omega_{m}}},\,
\displaystyle\overline{\omega_{m}}=\frac{d_m}{(IS_m+\epsilon)^5},
\end{equation*}
with linear weights $d_1=\dots=d_6=\displaystyle\frac{1}{6}$. The small positive constant is set to $\epsilon=10^{-15}$.
Details of the computation of $\chi$ can be found in \cite{GENO,CGKS5-prep}.

\subsection{Second-order gas-kinetic scheme}

For consistency, we retain the notation used in Section \ref{CGKS-5}. 
For GKS-2nd, the reconstruction stencil for the target cell $\Omega_0$ is defined as
\begin{equation*}
S^{cell-2nd}=\{\Omega_0, \Omega_{1,m}\,|\,m = 1, \dots, 6\},
\end{equation*}
which involves only face-neighboring cells, as shown in Figure \ref{stencil}(b).
The second-order reconstruction can be achieved based on
\begin{align*}
S^{data-2nd}=\{\boldsymbol{Q}_0,\,\boldsymbol{Q}_{1,m}|\,m = 1, \dots, 6\}.
\end{align*}
A linear polynomial $P^1(\boldsymbol{x})$ is constructed using $S^{data-2nd}$ and is defined as in Eq.\eqref{p1}.
To determine this second-order polynomial, the following system is solved using the least-square method
\begin{equation*}
\frac{1}{\left|\Omega_{k}\right|}\int_{\Omega_{k}}P^1(\boldsymbol{x})\text{d}V=\boldsymbol{Q}_{k},~\Omega_{k}\in S^{cell-2nd}.
\end{equation*}
To handle discontinuities, the linear polynomial $P^1(\boldsymbol{x})$ is coupled with the Venkatakrishnan limiter \cite{VK-limiter} to enable a nonlinear reconstruction as follows
\begin{equation*}
    \boldsymbol{Q}(\boldsymbol{x}_{G})=\boldsymbol{Q}_{0}+\phi_{0}\sum_{|\boldsymbol d|=1}b_{\boldsymbol d}p_{\boldsymbol d}(\boldsymbol{x}_{G}).
\end{equation*}
The limiter value $\phi_{0}$ for the target cell $\Omega_0$ is taken as
\begin{equation*}
    \phi_{0}=\min\{\phi_{F_m}\,|\,m = 1, \dots, 6\},
\end{equation*}
where $\phi_{F_m}$ denotes the Venkatakrishnan limiter \cite{VK-limiter} evaluated on face $F_m$.

\subsection{Temporal discretization}
    
According to Eq.\eqref{semi}, the discretization of the conservation laws over cell $\Omega_{i}$ can be represented as
\begin{equation*}
\boldsymbol{Q}_{i}^{n+1}=\boldsymbol{Q}_{i}^{n}+\int_{t^n}^{t^{n+1}} \mathcal{L}(\boldsymbol{Q}_{i},t)\,\mathrm{d}t,
\end{equation*}
where $\boldsymbol{Q}_{i}^{n+1}$ represents the cell-averaged conservative variable over cell $\Omega_{i}$ at $t_{n+1}=t_n+\Delta t$.
For GKS-2nd, a second-order temporal discretization is used
\begin{align*}
\boldsymbol{Q}_{i}^{n+1}=\boldsymbol{Q}_{i}^n+\Delta t\mathcal{L}(\boldsymbol{Q}_{i}^n)+\frac{1}{2}\Delta t^2\frac{\partial}{\partial t}\mathcal{L}(\boldsymbol{Q}^n_{i}).
\end{align*}
For the CGKS-5th scheme, high-order temporal accuracy is achieved with a two-stage, fourth-order method \cite{s2o4-0,GRP-high-1,GKS-high-1}
\begin{align*}
\boldsymbol{Q}_{i}^{*}&=\boldsymbol{Q}_{i}^n+\frac{1}{2}\Delta t\mathcal{L}(\boldsymbol{Q}_{i}^n)+\frac{1}{8}\Delta t^2\frac{\partial}{\partial t}\mathcal{L}(\boldsymbol{Q}^n_{i}), \\
\boldsymbol{Q}_{i}^{n+1}&=\boldsymbol{Q}_{i}^n+\Delta t\mathcal {L}(\boldsymbol{Q}_{i}^n)+\frac{1}{6}\Delta t^2\big(\frac{\partial}{\partial t}\mathcal{L}(\boldsymbol{Q}^n_{i})+2\frac{\partial}{\partial t}\mathcal{L}(\boldsymbol{Q}^{*}_{i})\big).
\end{align*}
The cell-averaged gradients and line-averaged derivatives are evaluated simultaneously using conservative variables at cell interfaces.
The evolution of the conservative variable on one side of an interface follows a two-stage update
\begin{align*}
\boldsymbol{Q}^{*}&=\boldsymbol{Q}^n+\frac{1}{2}\Delta t (\partial_t \boldsymbol{Q})^n,\\
\boldsymbol{Q}^{n+1}&=\boldsymbol{Q}^n+\Delta t (\partial_t \boldsymbol{Q})^{*}.
\end{align*}
To provide conservative variables on both sides of an interface, which may differ near discontinuities, the update model for $\boldsymbol{Q}(t)$ is
\begin{align*}
\boldsymbol{Q}^l(t)&=(1-e^{-\Delta t /\tau_0})\boldsymbol{Q}^e(t)+e^{-\Delta t /\tau_0}\boldsymbol{Q}_0^l(t),\\
\boldsymbol{Q}^r(t)&=(1-e^{-\Delta t /\tau_0})\boldsymbol{Q}^e(t)+e^{-\Delta t /\tau_0}\boldsymbol{Q}_0^r(t).
\end{align*}
Further details are provided in \cite{CGKS-high-5,CGKS-high-4}.

\section{Numerical verification}

In this section, the performance of CGKS-5th and GKS-2nd is evaluated using two approaches.
The first approach quantifies the computational resources required to achieve a prescribed target resolution.
The second approach evaluates the accuracy and resolution of the schemes under equivalent computational resources.
Using these criteria, numerical tests are conducted ranging from subsonic to supersonic turbulent flows, thereby assessing performance across smooth and discontinuous solutions. 
For clarity, CGKS-5th and GKS-2nd with linear reconstructions are specifically referred to as linear CGKS-5th and linear GKS-2nd, respectively.

For evolution of flow fields, the time step $\Delta t$ is given by the CFL condition.
For viscous flows, the time step is additionally constrained by the viscous term as $\Delta t = \text{CFL}\cdot h^2/(3\nu)$, where $h$ is the cell size and $\nu$
is the kinematic viscosity coefficient. 
For the inviscid flows, the collision time $\tau$ takes
\begin{align*}
\tau=\epsilon \Delta t+C\displaystyle|\frac{p_l-p_r}{p_l+p_r}|\Delta t,
\end{align*}
where $p_l$ and $p_r$ denote the pressure on the left and right sides of the cell interface, $\epsilon=0.05$ and $C=10.0$. For the viscous flows, the collision time is related to the viscosity coefficient,
\begin{align*}
\tau=\frac{\mu}{p}+C \displaystyle|\frac{p_l-p_r}{p_l+p_r}|\Delta t,
\end{align*}
where $\mu$ is the dynamic viscous coefficient and $p$ is the pressure at the cell interface. 
In smooth flow regions, the collision time reduces to
\begin{equation*}
\tau=\frac{\mu}{p}.
\end{equation*}

In the numerical cases presented in this section, CGKS-5th and GKS-2nd are implemented on multiple GPUs using the Compute Unified Device Architecture (CUDA) for parallel computations and the Message Passing Interface (MPI) for inter-process communication. This configuration enables large-scale simulations and is executed on NVIDIA GeForce RTX 4090 GPUs with double-precision.

\begin{figure}[!h]
\centering
\includegraphics[width=0.485\textwidth]{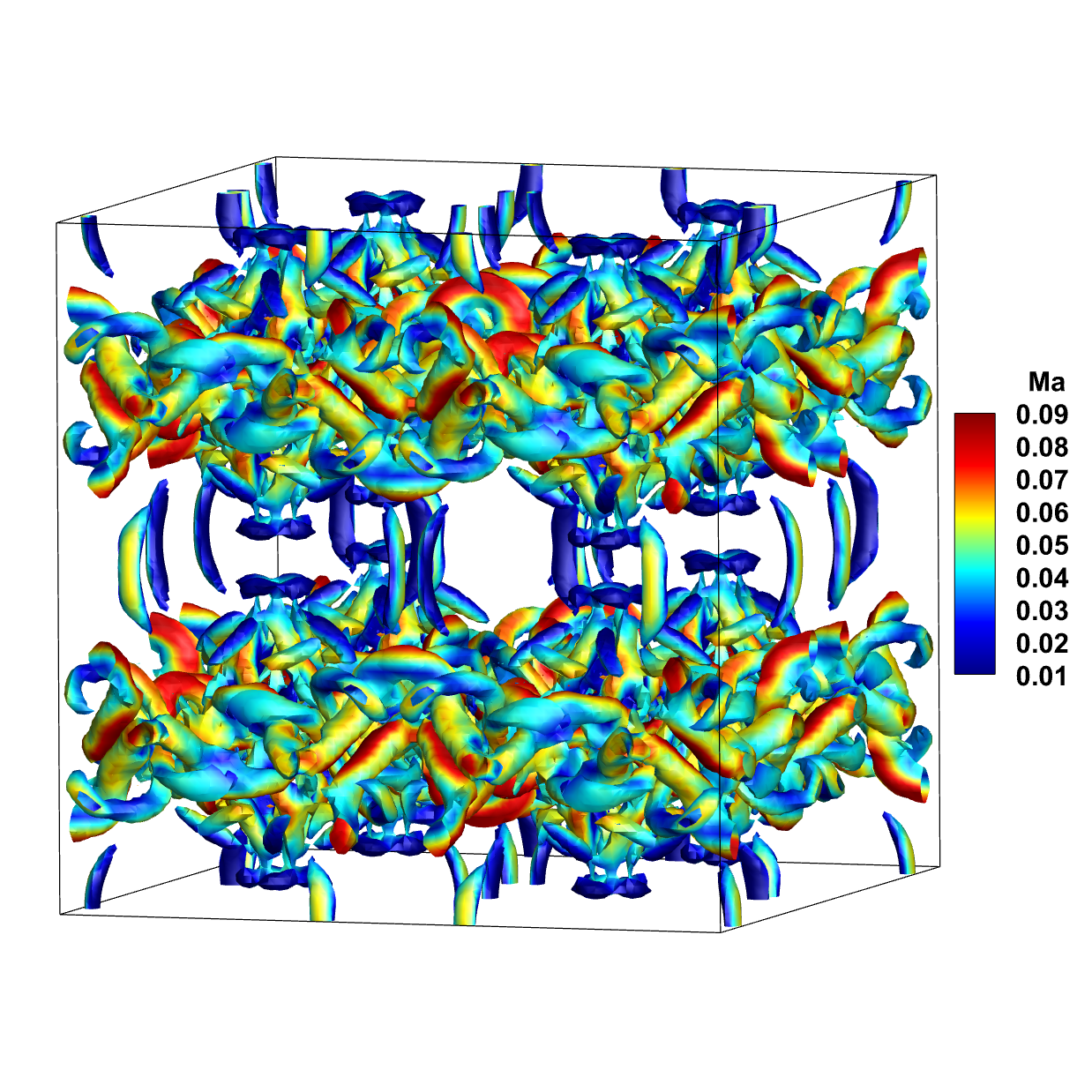}
\includegraphics[width=0.485\textwidth]{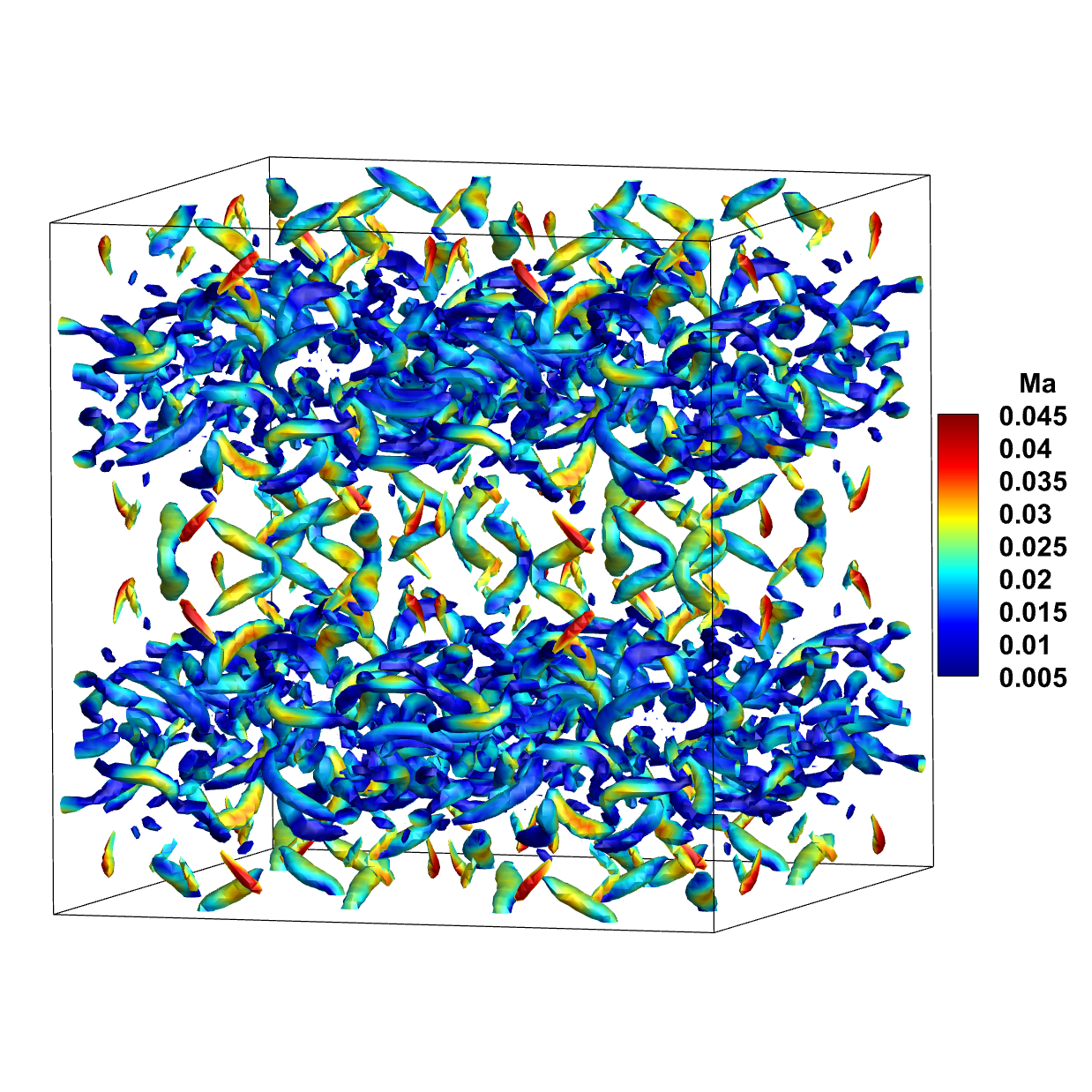}
\caption{\label{TGV-sub-1} Subsonic Taylor-Green vortex flow: Iso-surfaces of the Q-criterion ($Q=1.6$) at $t=10$ (left) and $t=20$ (right) for $Ma_{\infty}=0.1$ and $Re=1600$ using the linear CGKS-5th on a $96^3$ mesh.}
\end{figure}

\begin{figure}[!h]
\centering
\includegraphics[width=0.495\textwidth]{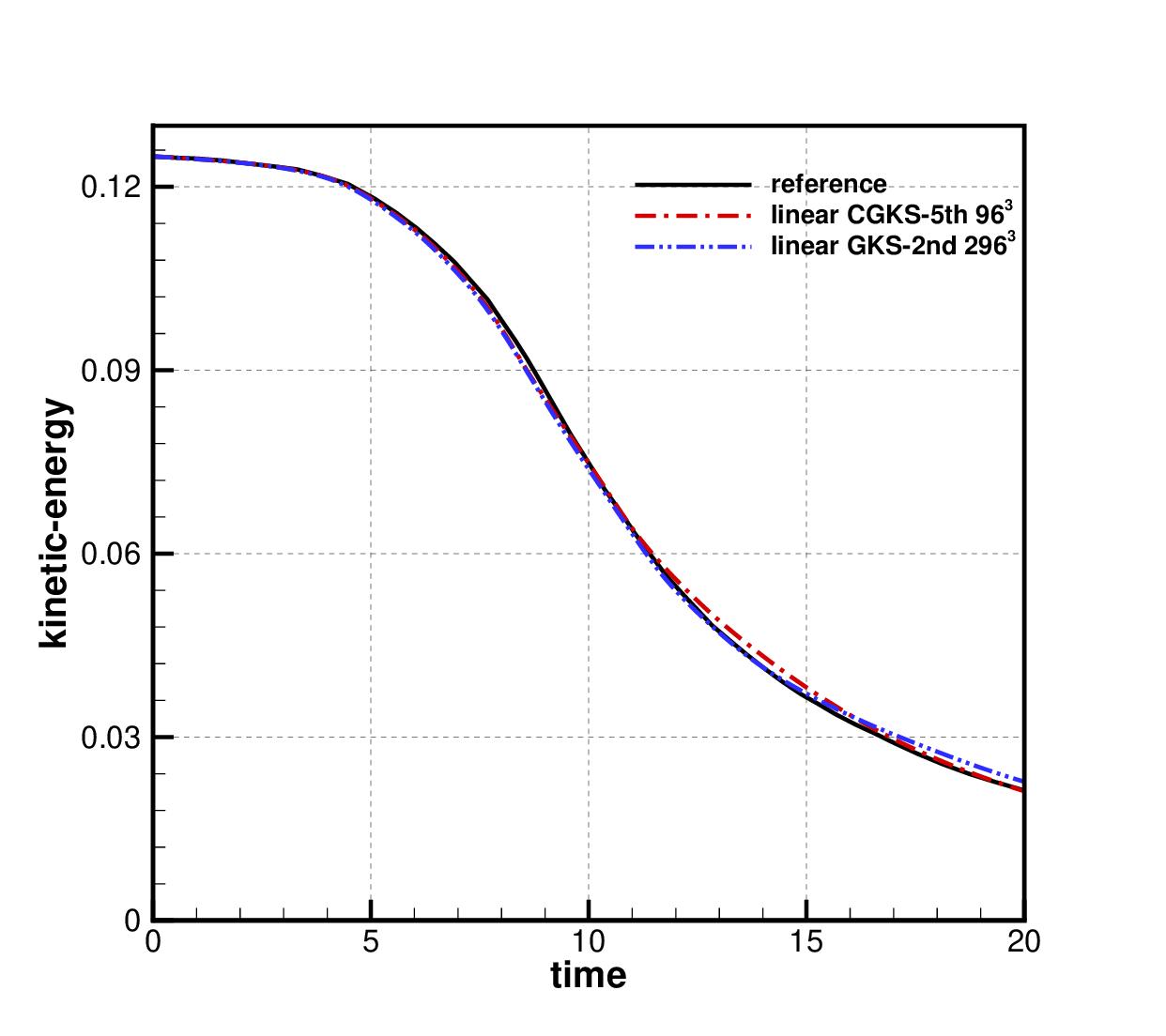}
\includegraphics[width=0.495\textwidth]{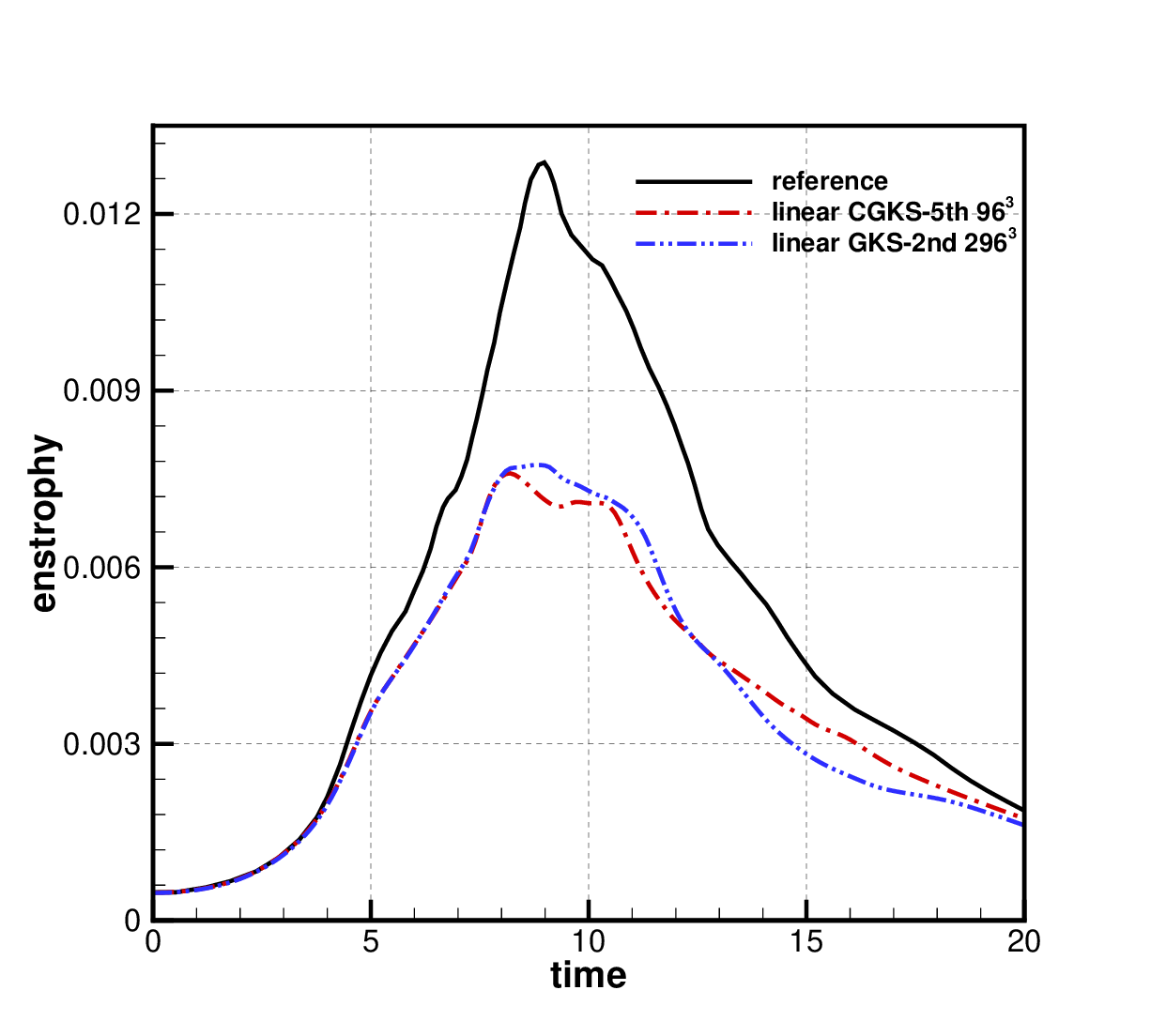}
\caption{\label{TGV-sub-2} Subsonic Taylor-Green vortex flow: the comparison of the time histories of kinetic energy and enstrophy dissipation with $Ma_{\infty}=0.1$ and $Re=1600$.}
\end{figure}

\subsection{Subsonic Taylor-Green vortex flow}

The Taylor-Green vortex is a well-established benchmark in fluid mechanics, widely used to assess numerical methods in terms of turbulence generation, energy dissipation, and the evolution of vortex structures \cite{TGV-1,TGV-2,TGV-3}.
In this subsection, the subsonic Taylor-Green vortex is selected as a representative benchmark case, enabling a detailed performance evaluation of the linear CGKS-5th and linear GKS-2nd. In this case, the comparison focuses on analyzing the computational cost required by each scheme to achieve a prescribed target resolution.

The initial velocity field is given by
\begin{align*}
U &= U_0\sin(\frac{x}{L})\cos(\frac{y}{L})\cos(\frac{z}{L}), \\
V &=-U_0\cos(\frac{x}{L})\sin(\frac{y}{L})\cos(\frac{z}{L}), \\
W &= 0, 
\end{align*}
and the initial pressure field is set as
\begin{equation}
\begin{aligned}\label{TGV-sub}
\rho &= \rho_0\big(1+\frac{\gamma Ma_{\infty}^2}{16}\big(\cos(\frac{2x}{L})+\cos(\frac{2y}{L})\big)\big(\cos(\frac{2z}{L})+2\big)\big),\\
p &= p_0\big(1+\frac{\gamma Ma_{\infty}^2}{16}\big(\cos(\frac{2x}{L})+\cos(\frac{2y}{L})\big)\big(\cos(\frac{2z}{L})+2\big)\big).
\end{aligned}
\end{equation}
In the computation, $U_0=1$, $\rho_0=1$.
The Reynolds number is $Re=1600$ defined by $Re=\rho_0U_0L/\mu_0$ , where $\mu_0$ is the dynamic viscosity coefficient. 
The Mach number is defined as $Ma_{\infty} = U_0 / \sqrt{\gamma p_0 / \rho_0}$, and this case is conducted at $Ma_{\infty} = 0.1$, which lies within the incompressible limit.
The Taylor-Green vortex flow is defined in a periodic domain defined as $[-\pi L, \pi L]^3$, $L=1$.
A uniform mesh with $96^3$ cells is utilized for the linear CGKS-5th, while a uniform mesh with $296^3$ cells is employed for the linear GKS-2nd.
It has been tested that the above grid configurations ensure comparable resolution for the two schemes.

Figure \ref{TGV-sub-1} presents the iso-surfaces of the Q-criterion ($Q=1.6$) at $t=10$ and $t=20$, computed using the linear CGKS-5th on a $96^3$ mesh, with the iso-surfaces colored by the Mach number. It can be observed that CGKS-5th is capable of capturing flow structures on a coarse mesh.
In order to quantitatively evaluate the performance of CGKS-5th and GKS-2nd, several time-averaged quantities are calculated.
The volume-averaged kinetic energy is given by
\begin{equation*}
E_k=\frac{1}{\rho_0 |\Omega|}\int_{\Omega} \frac{1}{2}\rho \boldsymbol{U}\cdot\boldsymbol{U}\,\mathrm{d}V,
\end{equation*}
where $\Omega$ is the computational domain and $E_k$ is calculated by numerical quadrature. The integrated enstrophy is defined as
\begin{equation*}
\epsilon^S=\frac{1}{\rho_0 |\Omega|}\int_{\Omega}\mu(|\nabla\times\boldsymbol{U}|^2)\,\mathrm{d}V,
\end{equation*}
where the dynamic viscosity $\mu$ is computed by the Sutherland's law, $\epsilon^S$ is calculated by numerical quadrature, and the velocity derivative values at quadrature points are obtained by the compact reconstruction.
The comparison of the time histories of kinetic energy $E_k$ and enstrophy $\epsilon^S$ using the linear CGKS-5th on a $96^3$ mesh and the linear GKS-2nd scheme on a $296^3$ mesh is presented in Figure \ref{TGV-sub-2}. The reference line in Figure \ref{TGV-sub-2} corresponds to a spectral solution computed on a $512^3$ mesh \cite{TGV-ref-1}.
As observed in Figure \ref{TGV-sub-2}, the two numerical results achieve equivalent resolution. 

Table \ref{TGV-sub-3} summarizes the computation time of the linear CGKS-5th on a $96^3$ mesh and the linear GKS-2nd on a $296^3$ mesh, both using multi-GPU parallel computing. As shown in Table \ref{TGV-sub-3}, CGKS-5th achieves nearly an order of magnitude improvement in computational efficiency compared to GKS-2nd for equivalent resolution. This demonstrates that, compared to the second-order GKS, the fifth-order CGKS can achieve higher resolution with significantly lower computational costs.

\begin{table}[!h]
\begin{center}
\def\temptablewidth{0.95\textwidth}{\rule{\temptablewidth}{1.0pt}}
\begin{tabular*}{\temptablewidth}{@{\extracolsep{\fill}}c|c|c|c|c} 
Number of GPUs & MESH  & Scheme & Computational time & Ratio  \\
\hline
2 & $96^3$ & linear CGKS-5th & 876 & 1 \\
2 & $296^3$ & linear GKS-2nd & 7746 & 8.8 \\
\end{tabular*}
{\rule{\temptablewidth}{1.0pt}}
\end{center}
\caption{\label{TGV-sub-3} Subsonic Taylor-Green vortex flow: the comparison of computational efficiency between linear CGKS-5th and linear GKS-2nd.}
\end{table}

\begin{figure}[!h]
\centering
\includegraphics[width=0.485\textwidth]{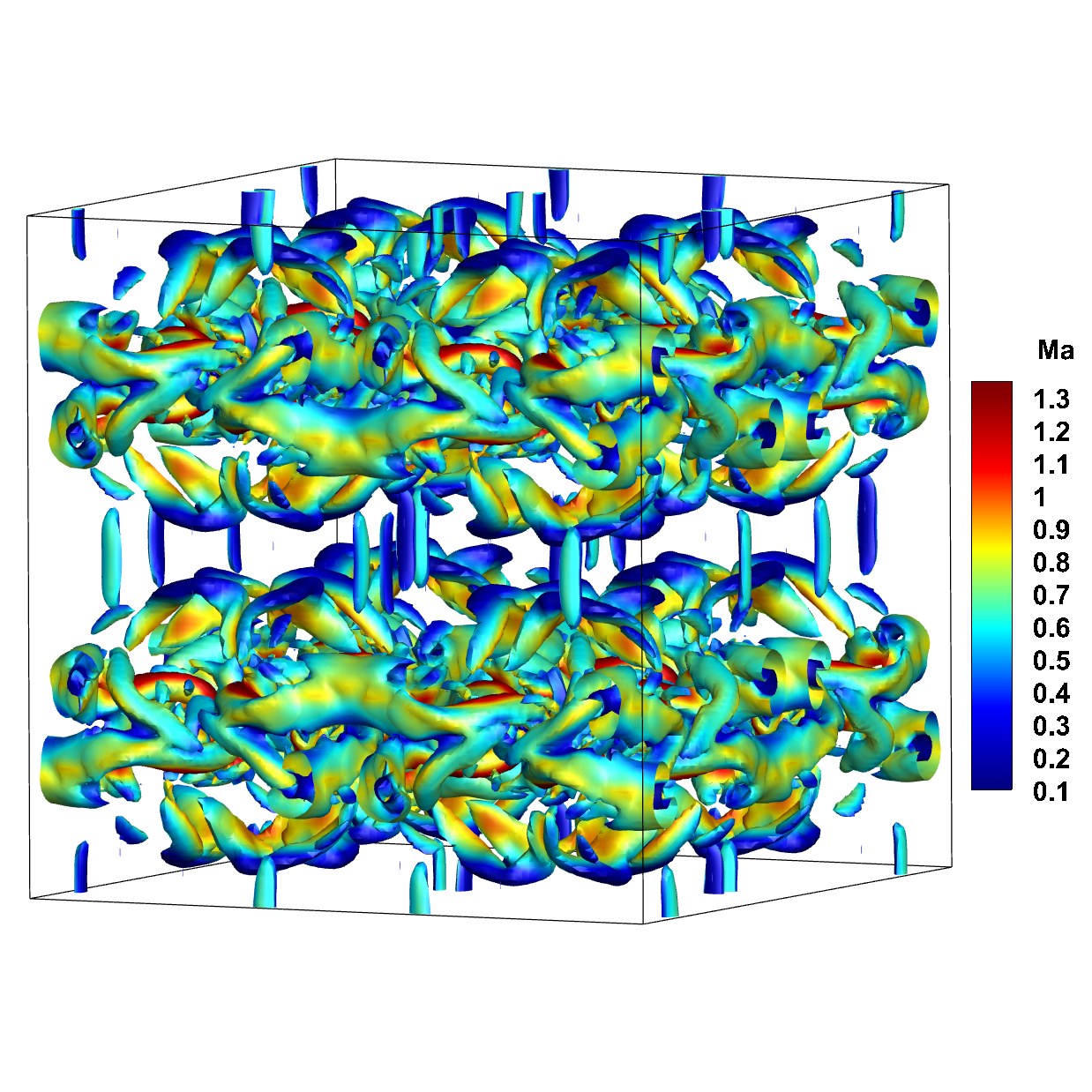}
\includegraphics[width=0.485\textwidth]{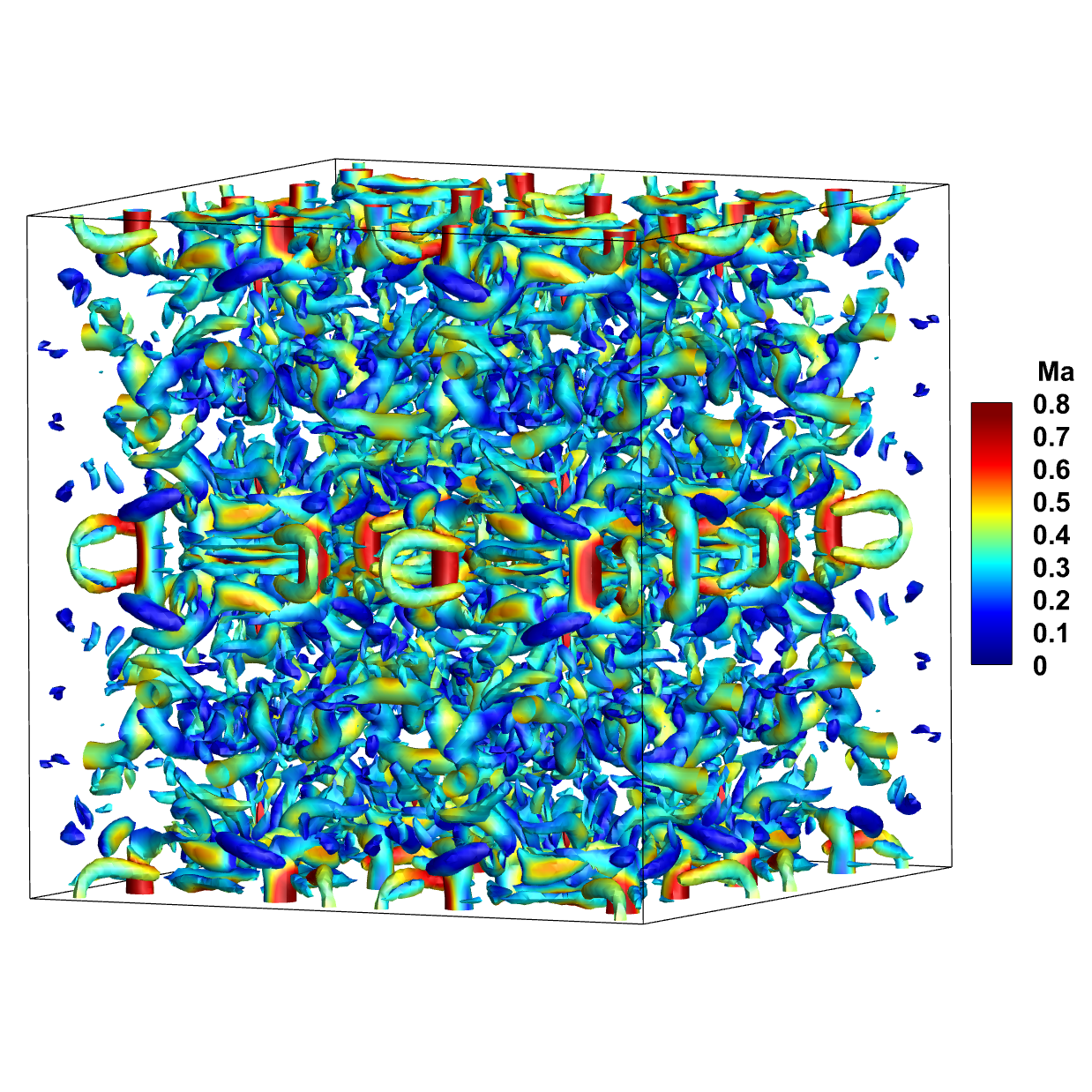}
\caption{\label{TGV-sup-1} High-speed Taylor-Green vortex flow: Iso-surfaces of the Q-criterion ($Q=4$) at $t=5$ (left) and $t=10$ (right) for $Ma_{\infty}=2.0$ and $Re=1600$ using CGKS-5th on a $116^3$ mesh.}
\includegraphics[width=0.48\textwidth]{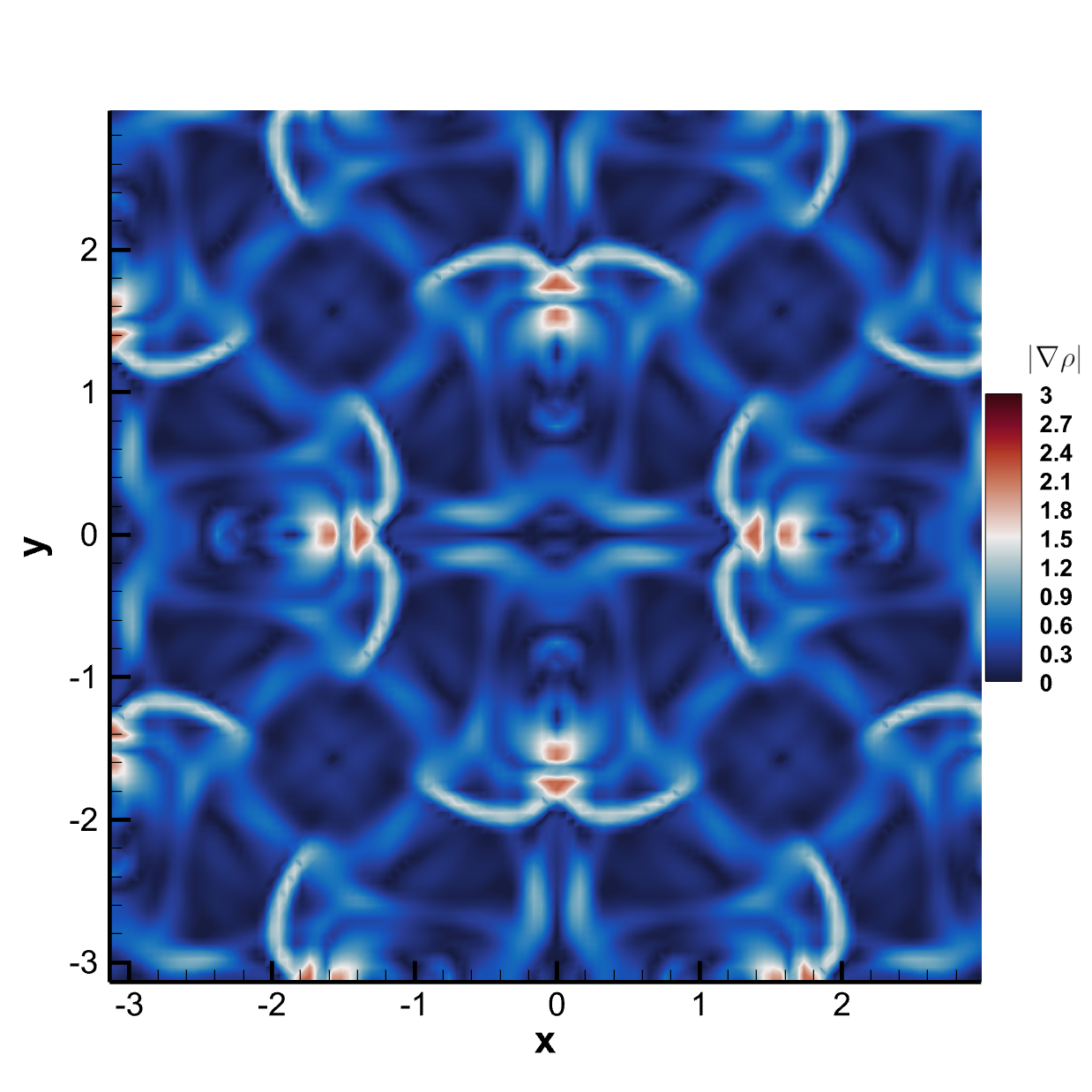}
\includegraphics[width=0.48\textwidth]{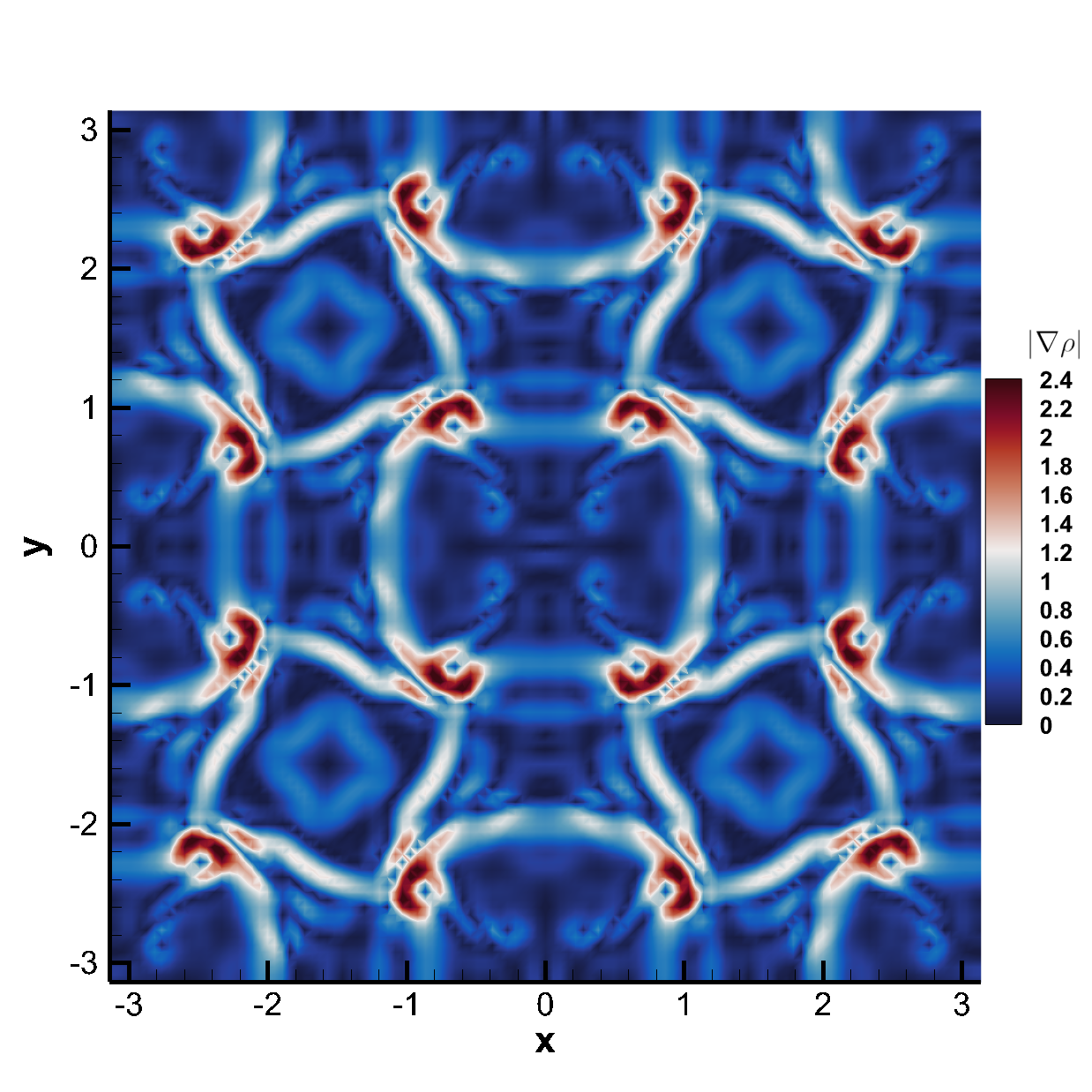}
\caption{\label{TGV-sup-2} High-speed Taylor-Green vortex flow: the magnitude of density gradients $\nabla\rho$ with $Ma_{\infty}=2.0$ and $Re=1600$ at $t=5$ (left) and $t=10$ (right) on the centerline (x-y) plane located at $z=0$ using CGKS-5th on a $116^3$ mesh.}
\end{figure}

\begin{figure}[!h]
\centering
\includegraphics[width=0.495\textwidth]{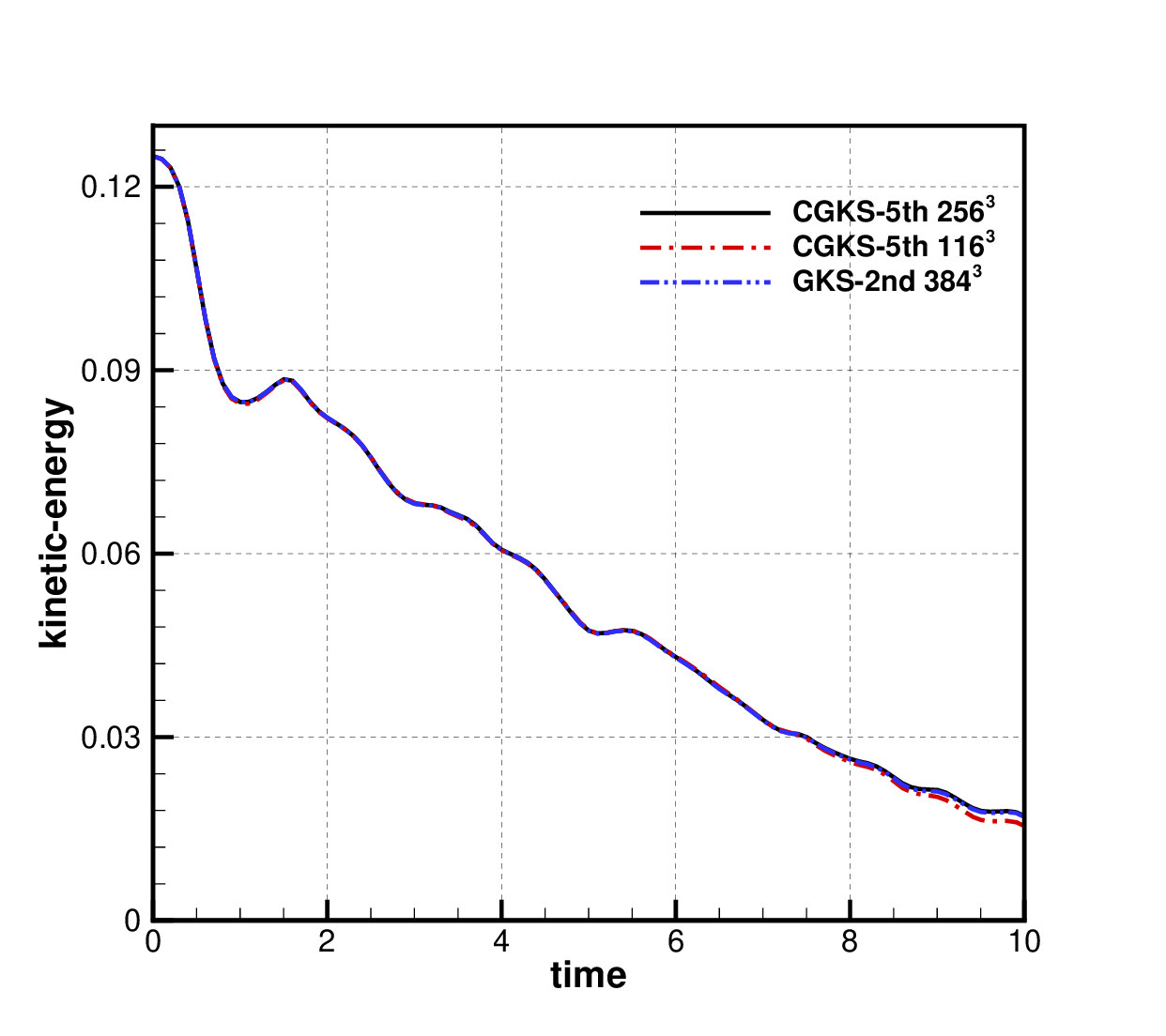}
\includegraphics[width=0.495\textwidth]{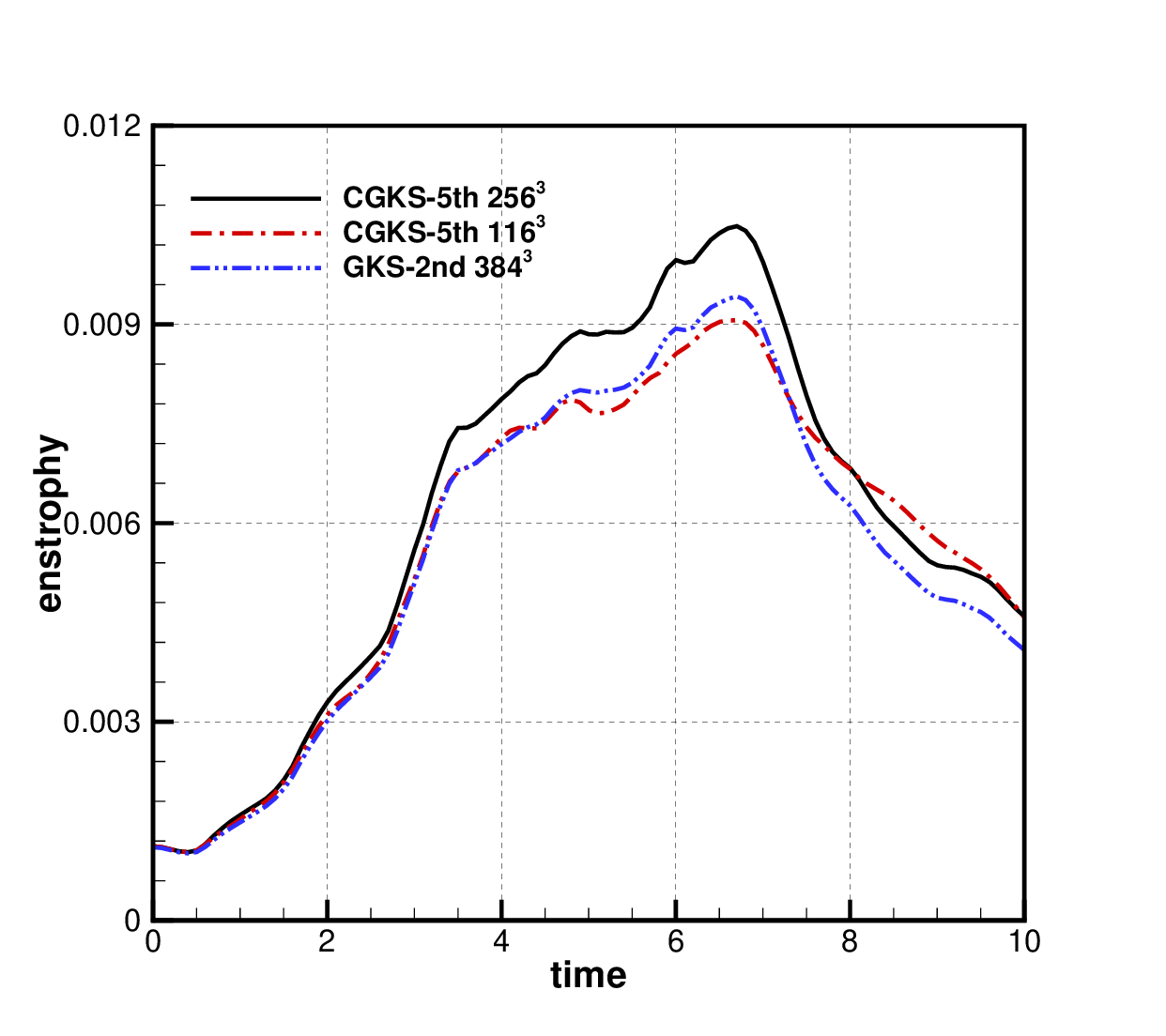}
\caption{\label{TGV-sup-3} High-speed Taylor-Green vortex flow: the comparison of the time histories of kinetic energy and enstrophy dissipation with $Ma_{\infty}=2.0$ and $Re=1600$.}
\end{figure}

\subsection{High-speed Taylor-Green vortex flow}

In previous studies, the subsonic and transonic Taylor-Green vortex flows have been extensively analyzed \cite{TGV-3, TGV-ref-4, TGV-ref-2}.
However, the high-speed Taylor-Green vortex flow (Mach number $Ma_{\infty} \geq 2$) has not been widely adopted as a benchmark problem.
Under such supersonic conditions, the original fluctuation form presented in Eq.\eqref{TGV-sub} cannot be directly extended to the high-speed regime.
To address this limitation, we propose a new initial condition, inspired by the classical form, that is tailored to high-speed flows.
In this subsection, the supersonic Taylor-Green vortex flow with $Ma_{\infty} = 2$ is considered, which exhibits pronounced shock-vortex interactions.
CGKS-5th and GKS-2nd employ nonlinear combination and Venkatakrishnan limiter to handle discontinuities, respectively, whereas the nonlinear part of CGKS-5th being algorithmically more complicated. The objective is to evaluate whether the nonlinear CGKS-5th can still deliver superior computational efficiency under these more demanding conditions.

The supersonic Taylor-Green vortex flow is defined in a periodic domain defined as $[-\pi L, \pi L]^3, L=1$. 
The initial condition is set as
\begin{align*}
U&=U_0\displaystyle\sin(\frac{x}{L})\cos(\frac{y}{L})\cos(\frac{z}{L}),\\
V&=-U_0\displaystyle\cos(\frac{x}{L})\sin(\frac{y}{L})\cos(\frac{z}{L}),\\
W&=0,
\end{align*}
and 
\begin{align*}
\rho&=\rho_0\displaystyle\big[1+\frac{1}{16}\big(\cos(\frac{2x}{L})+\cos(\frac{2y}{L})\big)\big(\cos(\frac{2z}{L})+2\big)\big],\\
p&=p_0\displaystyle\big[1+\frac{1}{16}\big(\cos(\frac{2x}{L})+\cos(\frac{2y}{L})\big)\big(\cos(\frac{2z}{L})+2\big)\big],\\
T&=1,\, R=1.
\end{align*}
In the computation, $\rho_0=1,\, p_0=1,\, U_0=\sqrt\gamma Ma_{\infty},\, Ma_{\infty}=2.0,\, T_0=1$, the Reynolds number $Re=1600$, and the Prandtl number is set as 0.7.
The dynamic viscosity $\mu$ is computed through the Sutherland law
\begin{equation*}
\mu(T)=\mu_0\displaystyle\frac{1.4042T^{3/2}}{T+0.4042},
\end{equation*}
where
$$\mu_0=\displaystyle\frac{\rho_0 U_0 L}{Re},\,T=\displaystyle\frac{p}{\rho}.$$
The kinetic energy integrated over the domain is given by
\begin{equation*}
E_k=\displaystyle\frac{1}{\rho_0 U_0^2 (2\pi)^3}\int_\Omega\frac{1}{2}\rho U\cdot U\, \mathrm{d}V,
\end{equation*}
where $\Omega$ is the computational domain and $E_k$ is calculated by numerical quadrature. 
The total viscous dissipation rate is defined as
\begin{align*}
\epsilon_T&=\epsilon_S+\epsilon_D,\\
&=\displaystyle\frac{1}{\rho_0 U_0^2 (2\pi)^3}\int_\Omega \mu(\nabla\times U)^2\, \mathrm{d}V+\displaystyle\frac{4}{3\rho_0 U_0^2 (2\pi)^3}\int_\Omega \mu(\nabla\cdot U)^2\, \mathrm{d}V.
\end{align*}
where $\epsilon^S$ is the solenoidal (enstrophy) dissipation, and $\epsilon^D$ is the dilatational dissipation, $\epsilon_T, \epsilon^S, \epsilon_D$ are calculated by numerical quadrature, and the velocity derivative values at quadrature points are obtained by the compact reconstruction.

A uniform mesh with $116^3$ cells is utilized for CGKS-5th, while a uniform mesh with $384^3$ cells is employed for GKS-2nd.
Figure \ref{TGV-sup-1} presents the iso-surfaces of the Q-criterion ($Q=4$) at $t=5$ and $t=10$, obtained using CGKS-5th on the $116^3$ mesh.
The iso-surfaces are colored by the Mach number, illustrating the flow structures.
Figure \ref{TGV-sup-2} shows the density gradient magnitude on the centerline (x-y) plane at $z=0$ for $t=5$ and $t=10$, computed using CGKS-5th on a $116^3$ mesh.
These results demonstrate the vortex decay over time and highlight the capability of CGKS-5th to accurately resolve flow structures with high resolution.
The comparison of the time histories of kinetic energy $E_k$ and enstrophy dissipation $\epsilon^S$ obtained using CGKS-5th on a $116^3$ mesh and GKS-2nd on a $384^3$ mesh is presented in Figure \ref{TGV-sup-3}. Additionally, the results from CGKS-5th on a $256^3$ mesh are used as the reference solution.
As shown in Figure \ref{TGV-sup-3}, the numerical results computed with CGKS-5th on a $116^3$ mesh achieve the same resolution as those computed with GKS-2nd on a $384^3$ mesh.
Table \ref{TGV-sup-4} compares the computational time required for these simulations using multi-GPU parallelization.
Notably, despite the higher algorithmic complexity of CGKS-5th caused by its nonlinear reconstruction, it still demonstrates a significant efficiency advantage over the limiter-based GKS-2nd, achieving an improvement close to an order of magnitude.

\begin{table}
\begin{center}
\def\temptablewidth{0.95\textwidth}{\rule{\temptablewidth}{1.0pt}}
\begin{tabular*}{\temptablewidth}{@{\extracolsep{\fill}}c|c|c|c|c} 
Number of GPUs & MESH  & Scheme & Computational time  & Ratio  \\
\hline
4 & $116^3$ & nonlinear CGKS-5th & 407 & 1 \\
4 & $384^3$ & nonlinear GKS-2nd & 2842 & 7.0 \\
\end{tabular*}
{\rule{\temptablewidth}{1.0pt}}
\end{center}
\caption{\label{TGV-sup-4} High-speed Taylor-Green vortex flow: the comparison of computational efficiency between CGKS-5th and GKS-2nd.}
\end{table}

\begin{figure}[!h]
\centering
\includegraphics[width=0.72\textwidth]{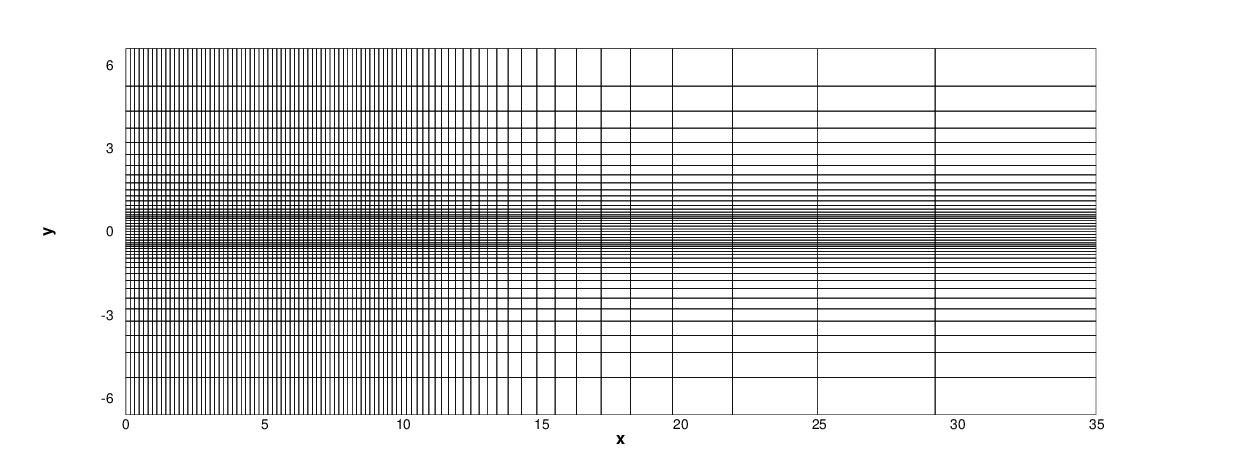}
\includegraphics[width=0.27\textwidth]{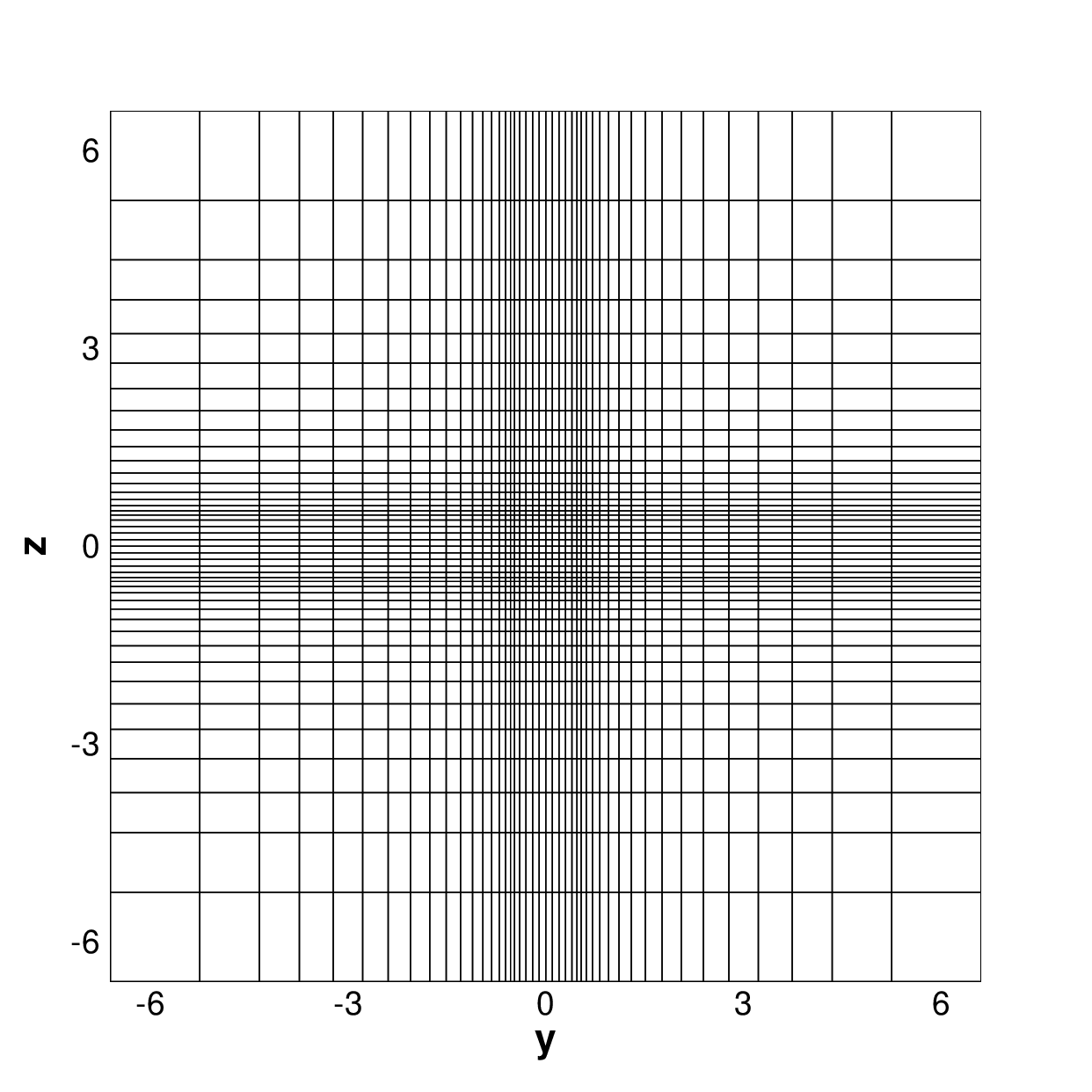}
\caption{\label{jet-mesh} Transonic Jet Simulation: Visualization of the computational mesh on the $y=0$ plane (left) and the $x=0$ plane (right), with every fourth line displayed in each of the three coordinate directions.}
\end{figure}

\begin{figure}[!h]
\centering
\includegraphics[width=0.7\textwidth]{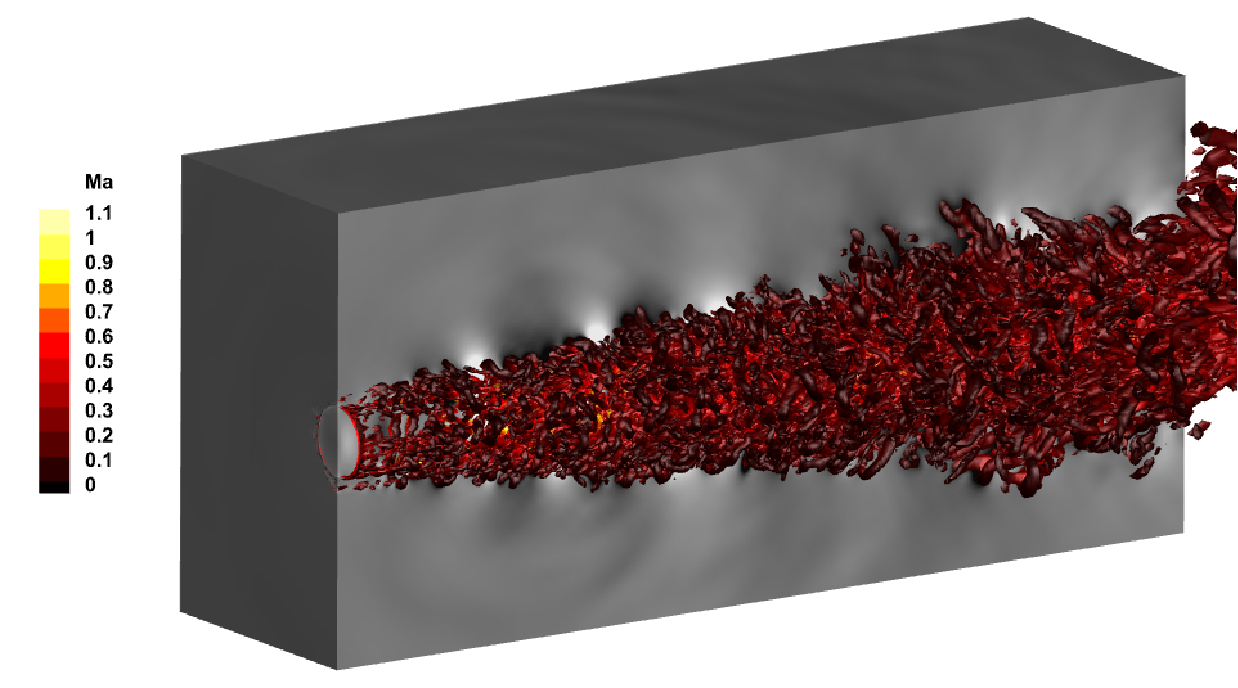}
\caption{\label{jet-sub-c5-Q} Transonic Jet Simulation: Isometric view of the iso-surface of the Q-criterion ($Q=0.2$) computed using CGKS-5th, accompanied by grayscale visualization of pressure fluctuations around $P_0 \pm 10^{-2}P_0$.}
\includegraphics[width=0.7\textwidth]{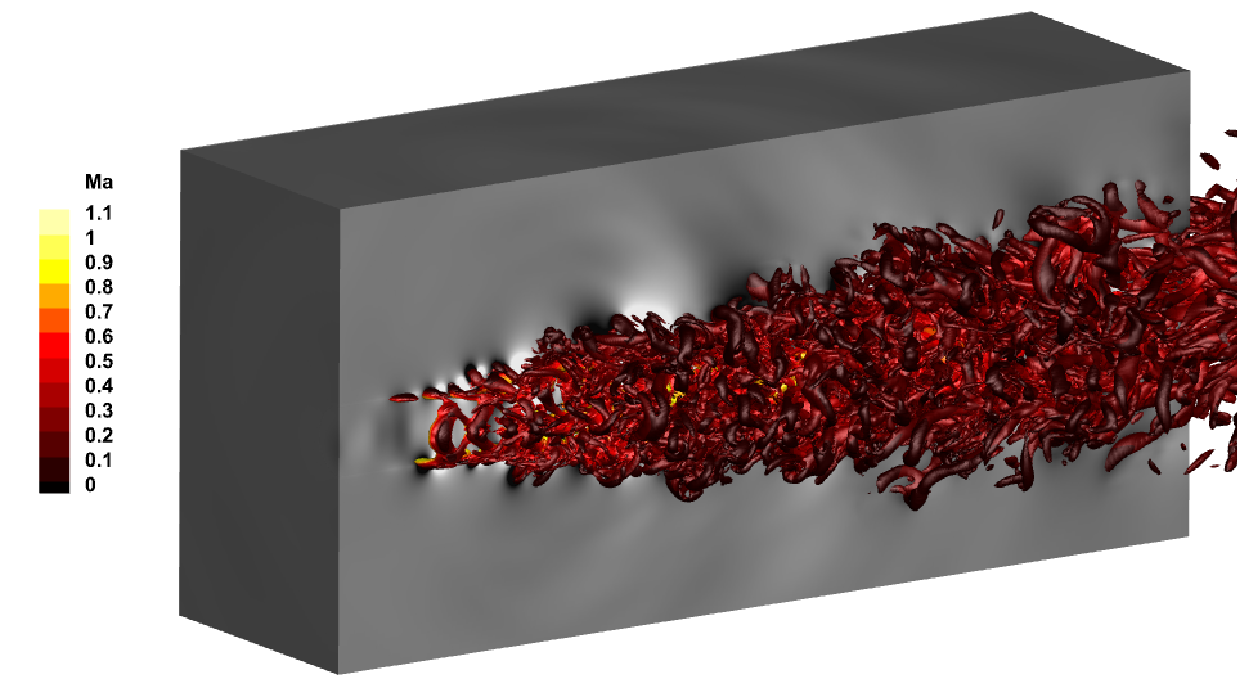}
\caption{\label{jet-sub-p1-Q} Transonic Jet Simulation: Isometric view of the iso-surface of the Q-criterion ($Q=0.2$) computed using GKS-2nd, accompanied by grayscale visualization of pressure fluctuations around $P_0 \pm 10^{-2}P_0$.}
\end{figure}

\begin{figure}[!h]
\centering
\includegraphics[width=0.485\textwidth]{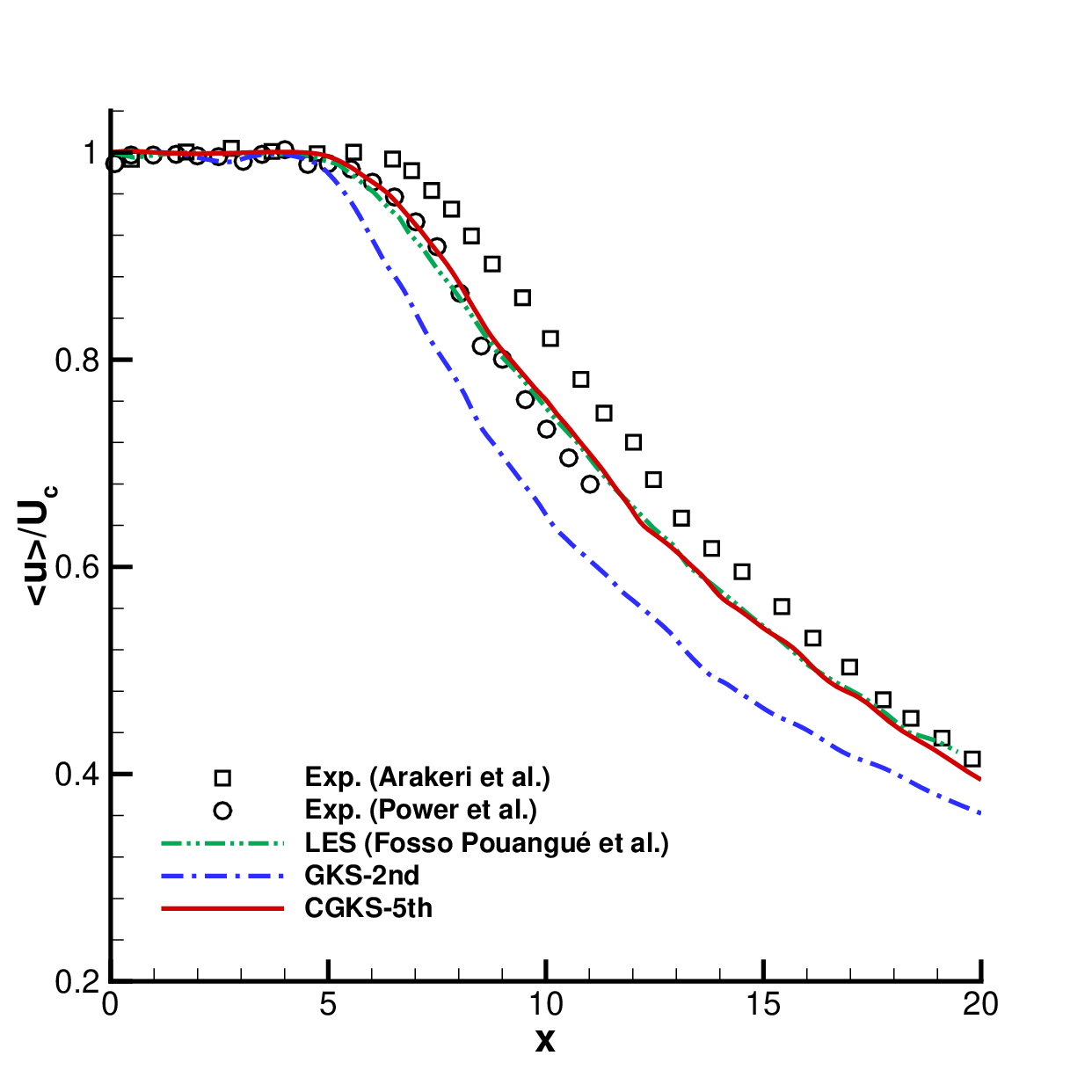}
\includegraphics[width=0.485\textwidth]{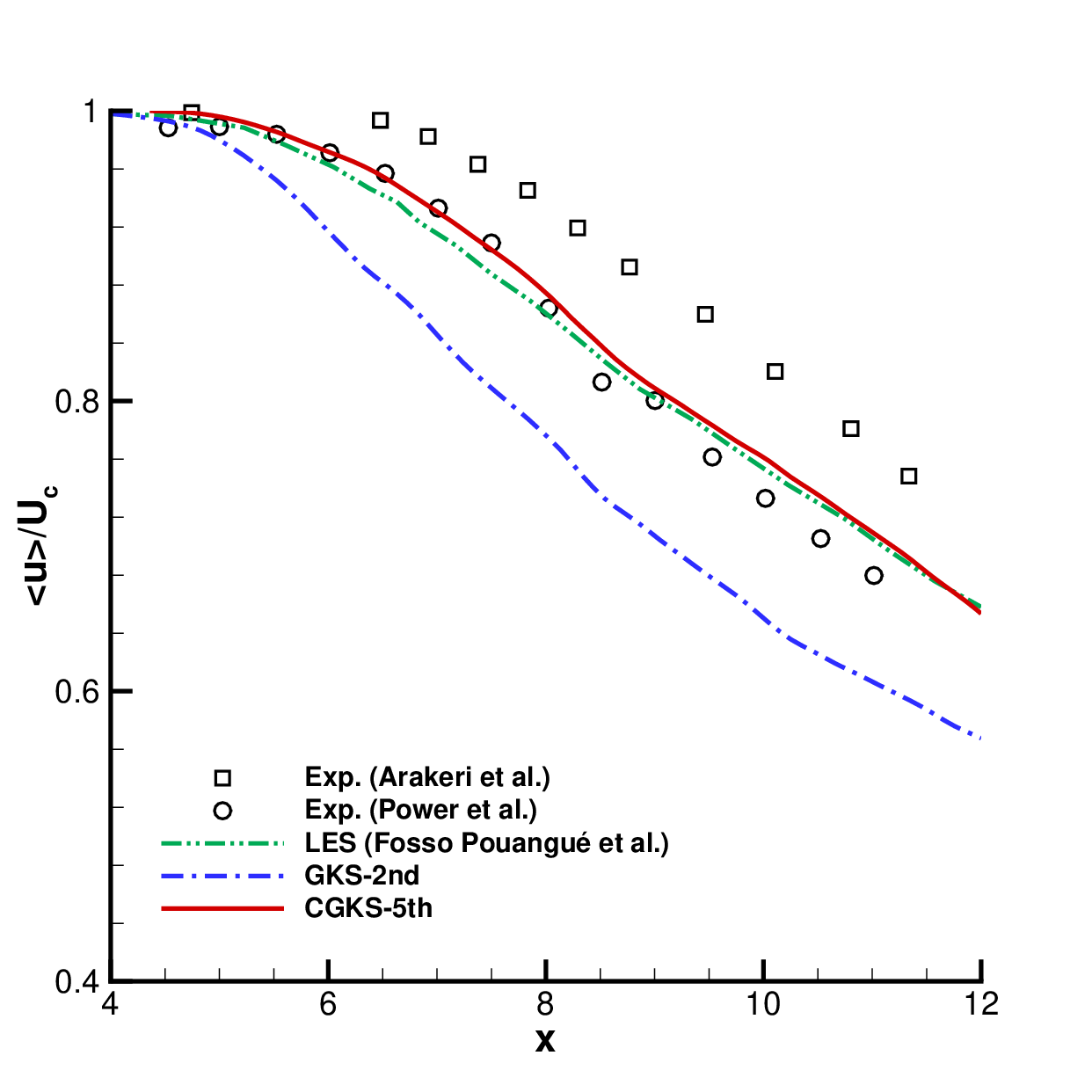}
\caption{\label{jet-sub-ave1} Transonic Jet Simulation: the mean axial velocity (left) and local magnification (right) along the centerline.}
\includegraphics[width=0.485\textwidth]{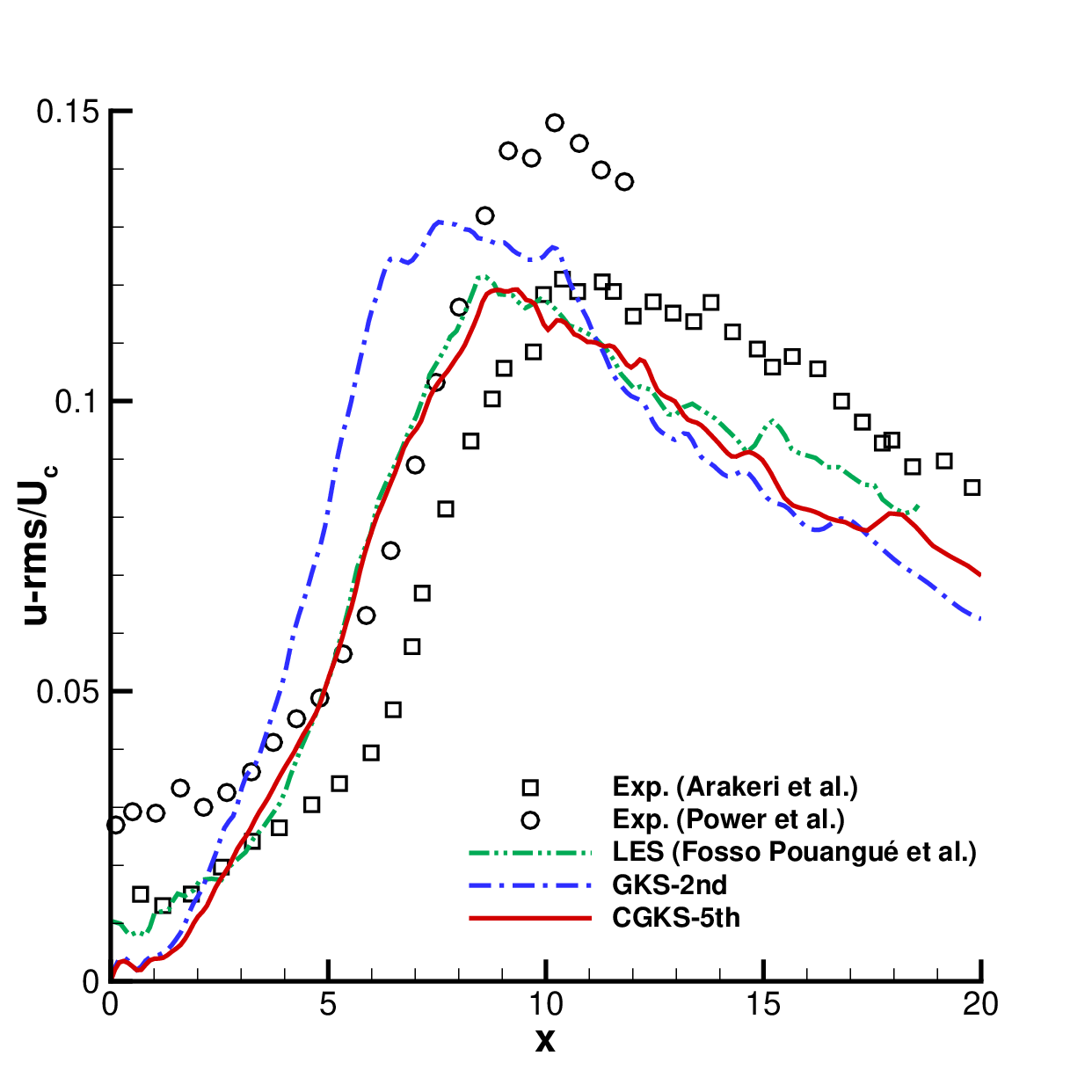}
\includegraphics[width=0.485\textwidth]{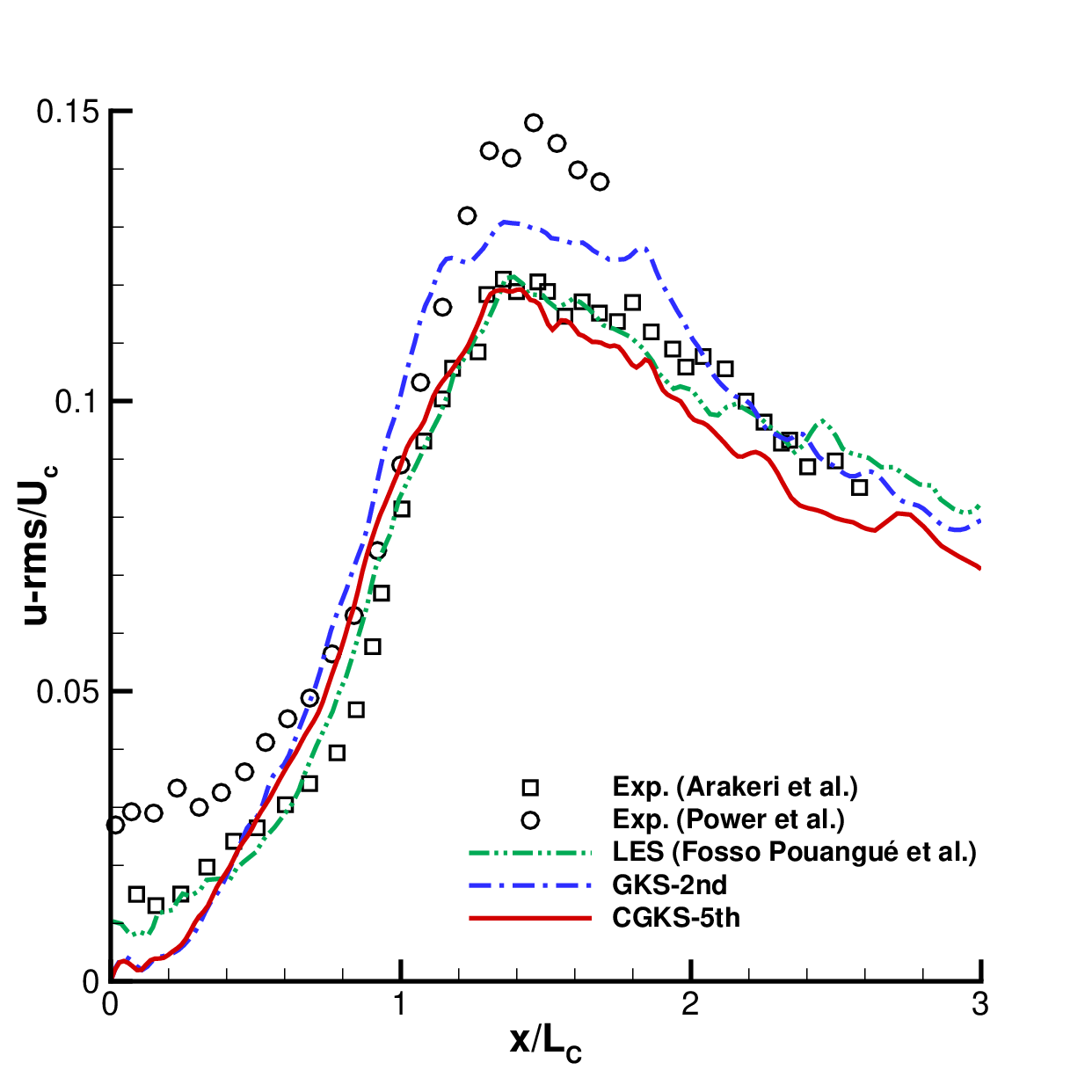}
\caption{\label{jet-sub-ave2} Transonic Jet Simulation: the root mean square axial velocity (left) and $L_C$-normalized distribution (right) along the centerline.}
\end{figure}

\begin{figure}[!h]
\centering
\includegraphics[width=0.7\textwidth]{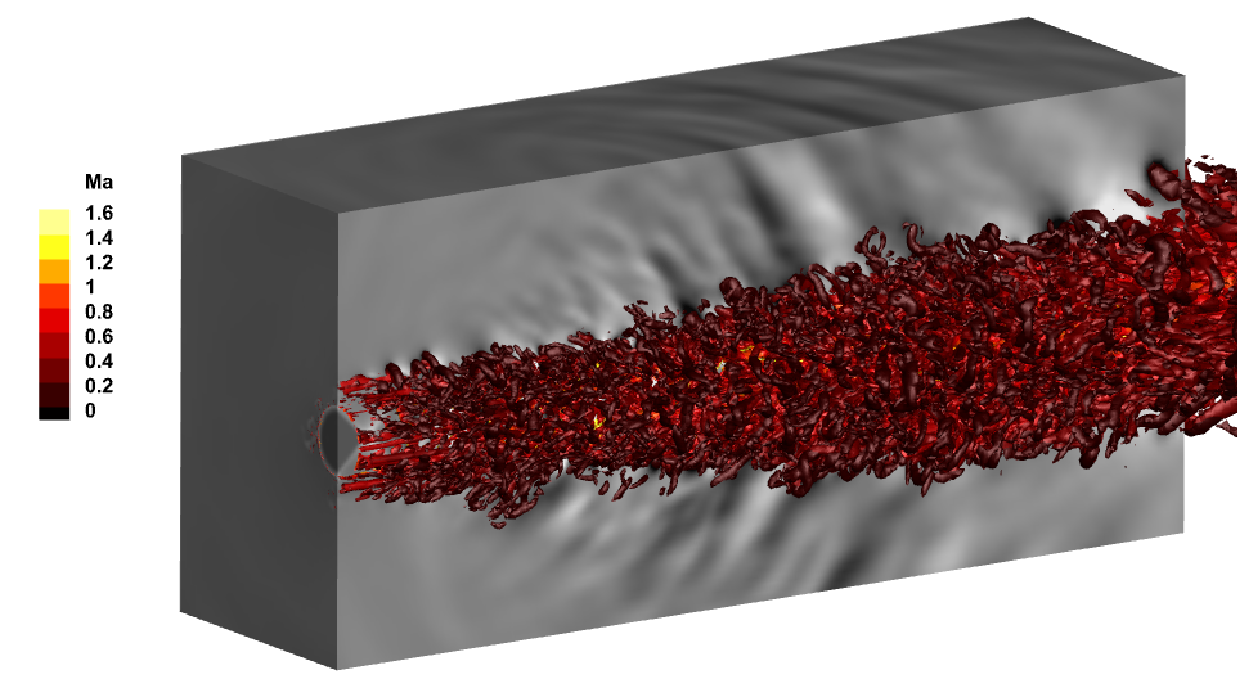}
\caption{\label{jet-sup-c5-Q} Supersonic Jet Simulation: the iso-surface of the Q-criterion ($Q=0.4$) using CGKS-5th shown in an isometric view, along with the pressure fluctuations around $P_0 \pm 4\times 10^{-2}P_0$ visualized in grayscale.}
\includegraphics[width=0.7\textwidth]{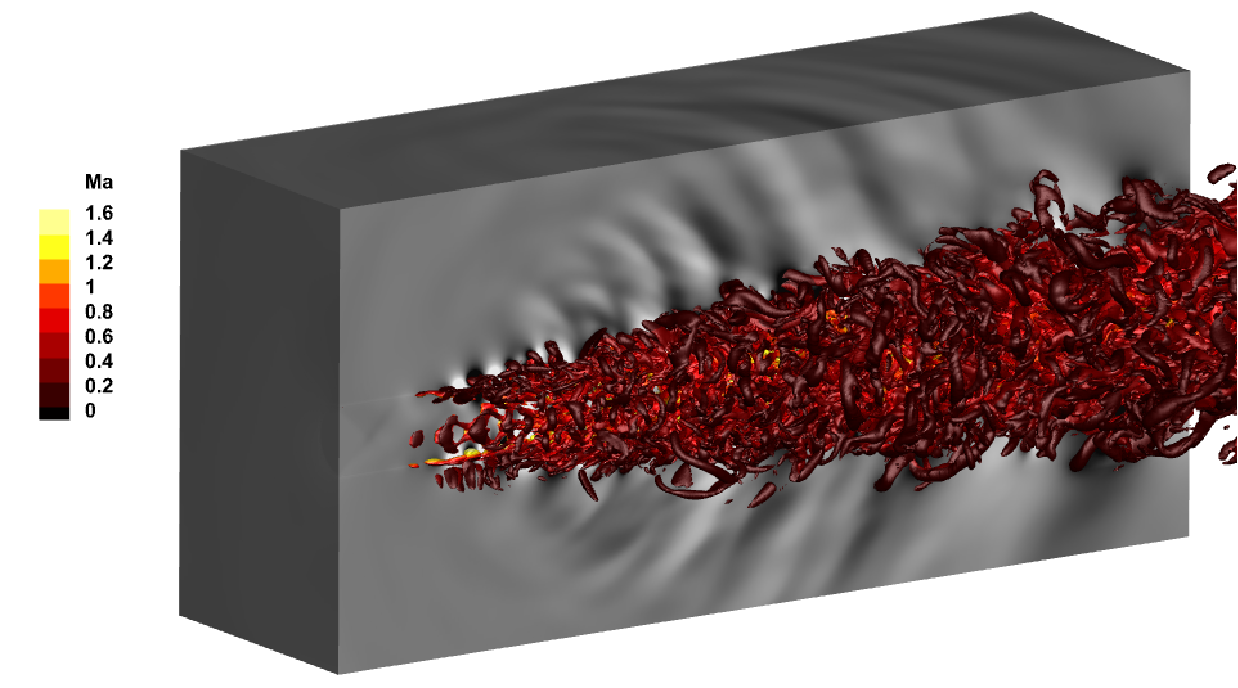}
\caption{\label{jet-sup-p1-Q} Supersonic Jet Simulation: the iso-surface of the Q-criterion ($Q=0.4$) using GKS-2nd shown in an isometric view, along with the pressure fluctuations around $P_0 \pm 4\times 10^{-2}P_0$ visualized in grayscale.}
\end{figure}

\begin{figure}[!h]
\centering
\includegraphics[width=0.485\textwidth]{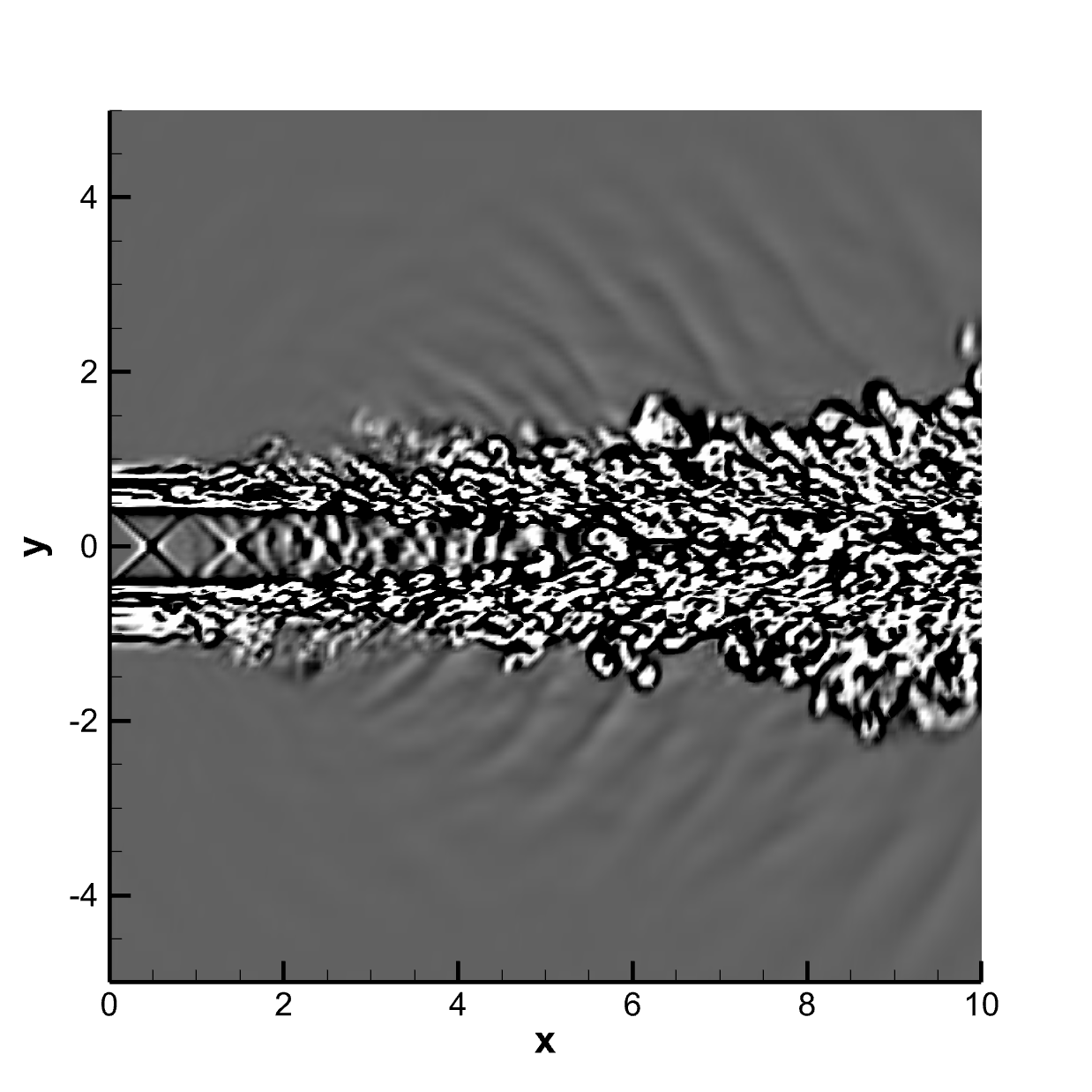}
\includegraphics[width=0.485\textwidth]{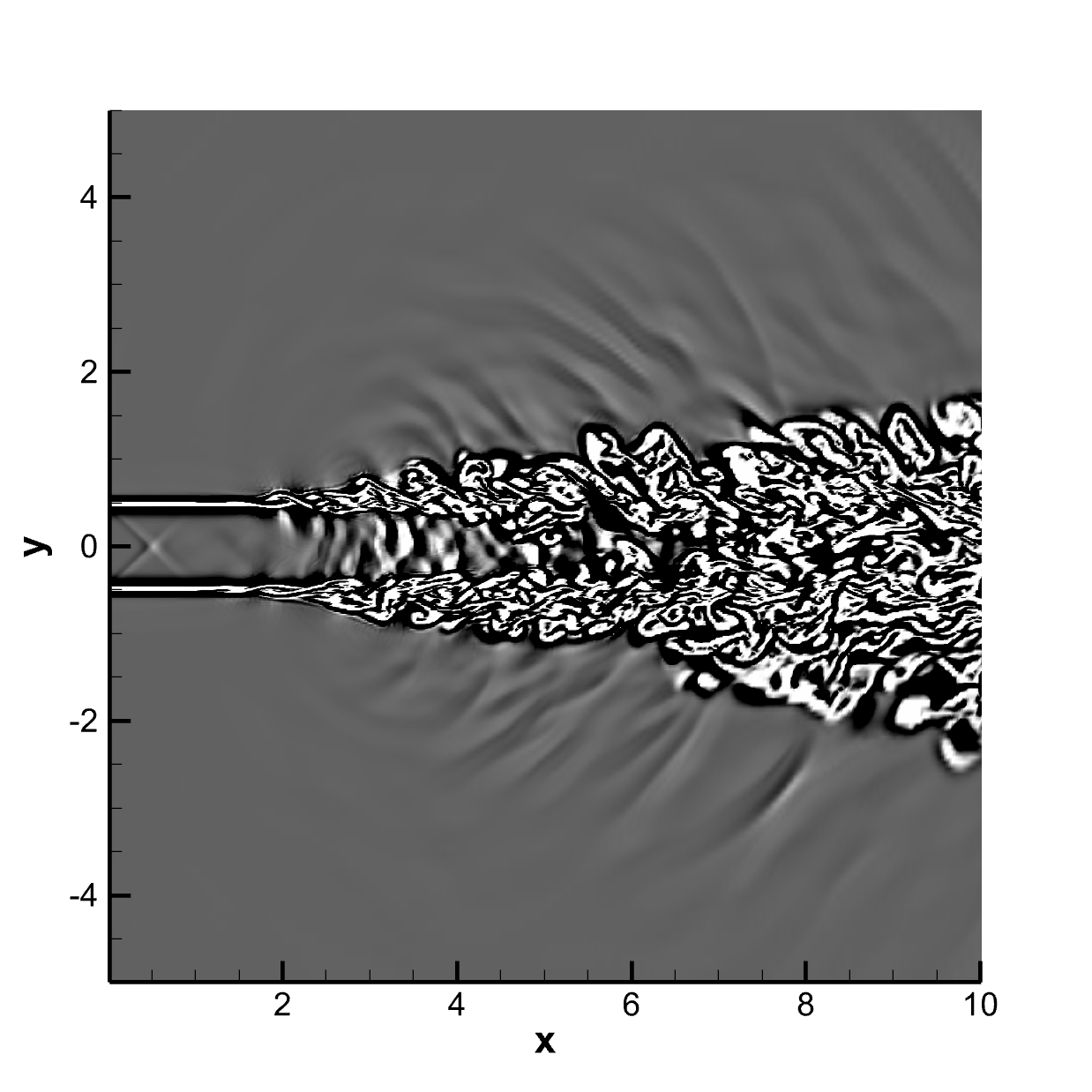}
\caption{\label{jet-sup-shadow} Supersonic Jet Simulation: the instantaneous shadowgraph using CGKS-5th (left) and GKS-2nd (right).}
\end{figure}

\begin{figure}[!h]
\centering
\includegraphics[width=0.485\textwidth]{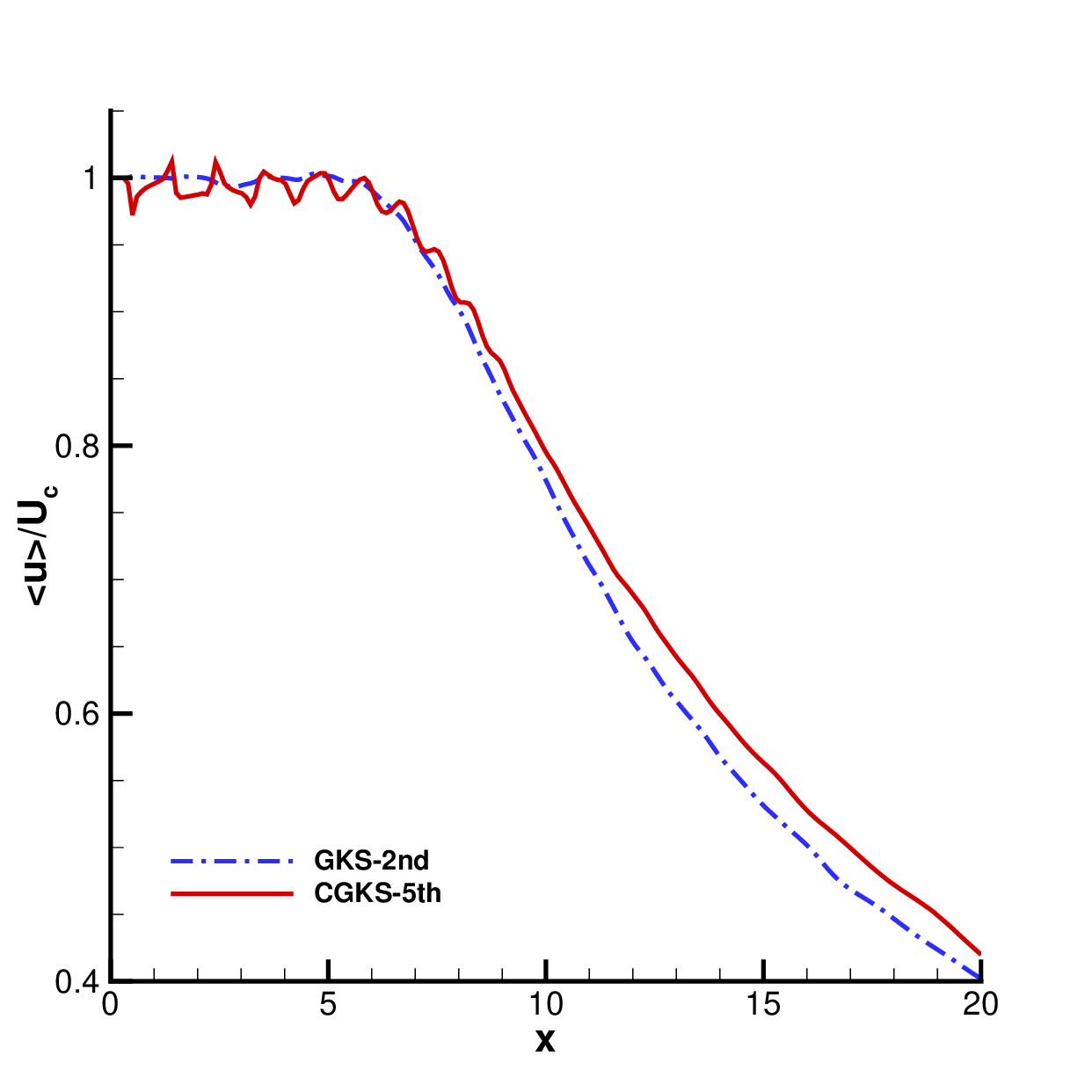}
\includegraphics[width=0.485\textwidth]{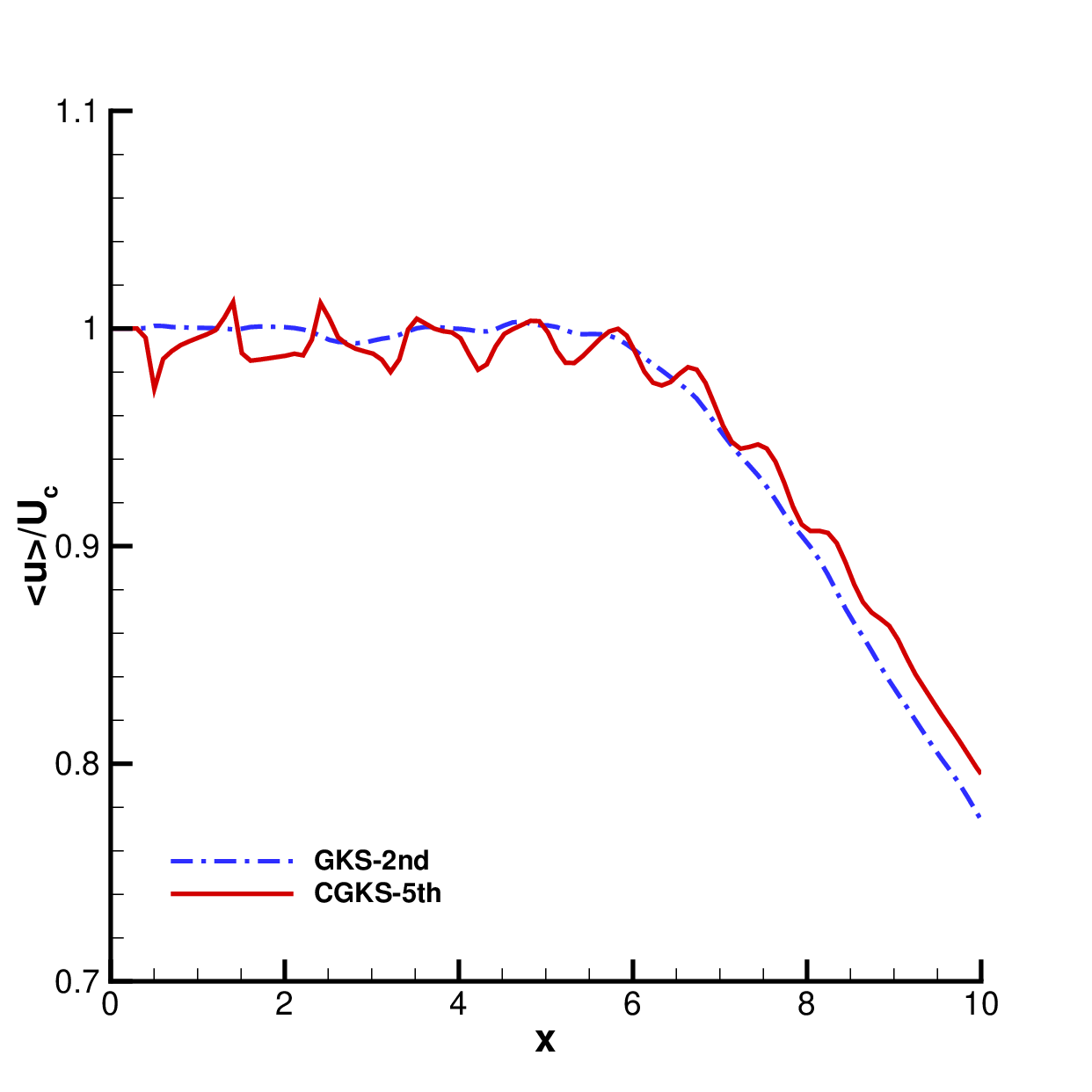}
\caption{\label{jet-sup-ave1} Supersonic Jet Simulation: the mean axial velocity (left) and local magnification (right) along the centerline.}
\includegraphics[width=0.485\textwidth]{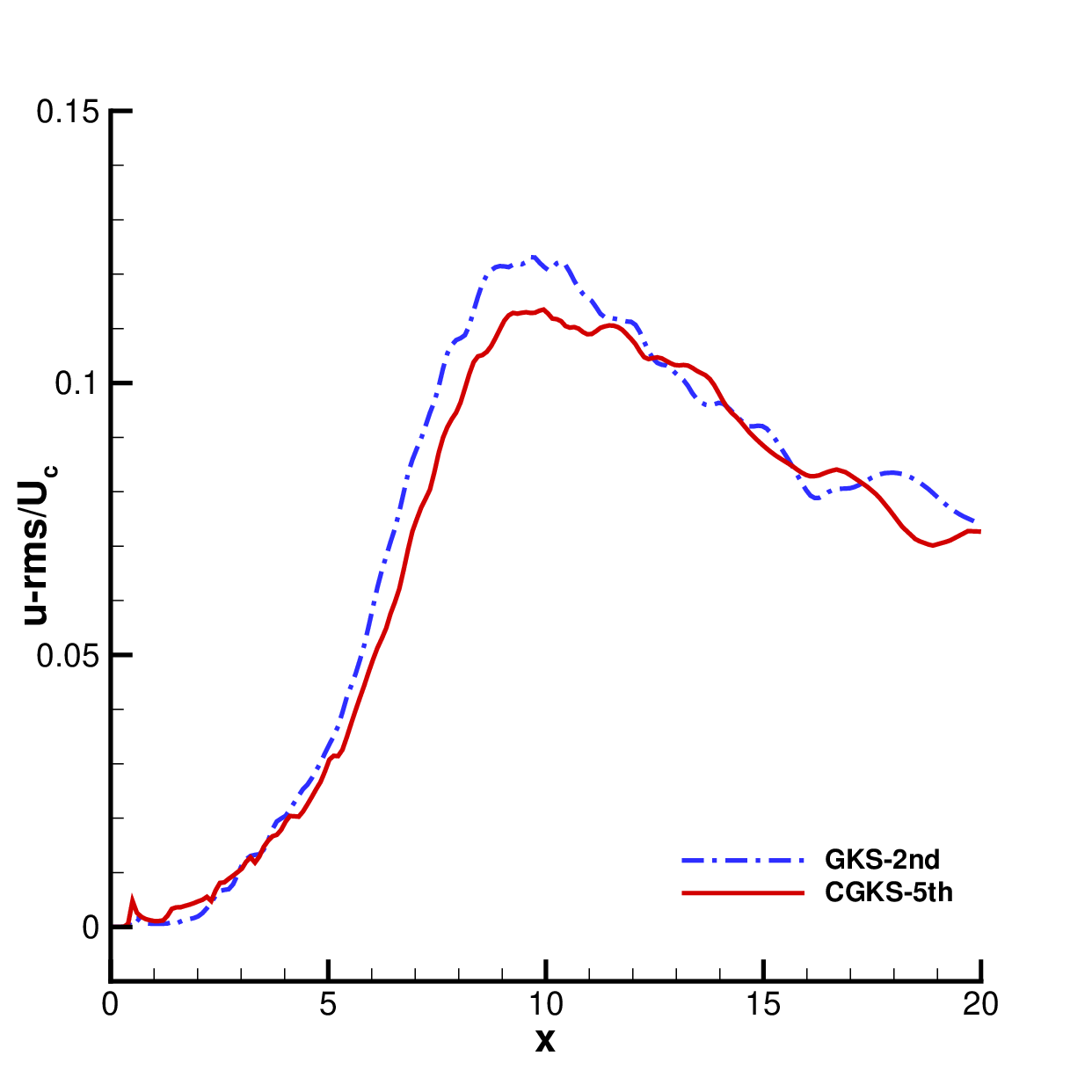}
\includegraphics[width=0.485\textwidth]{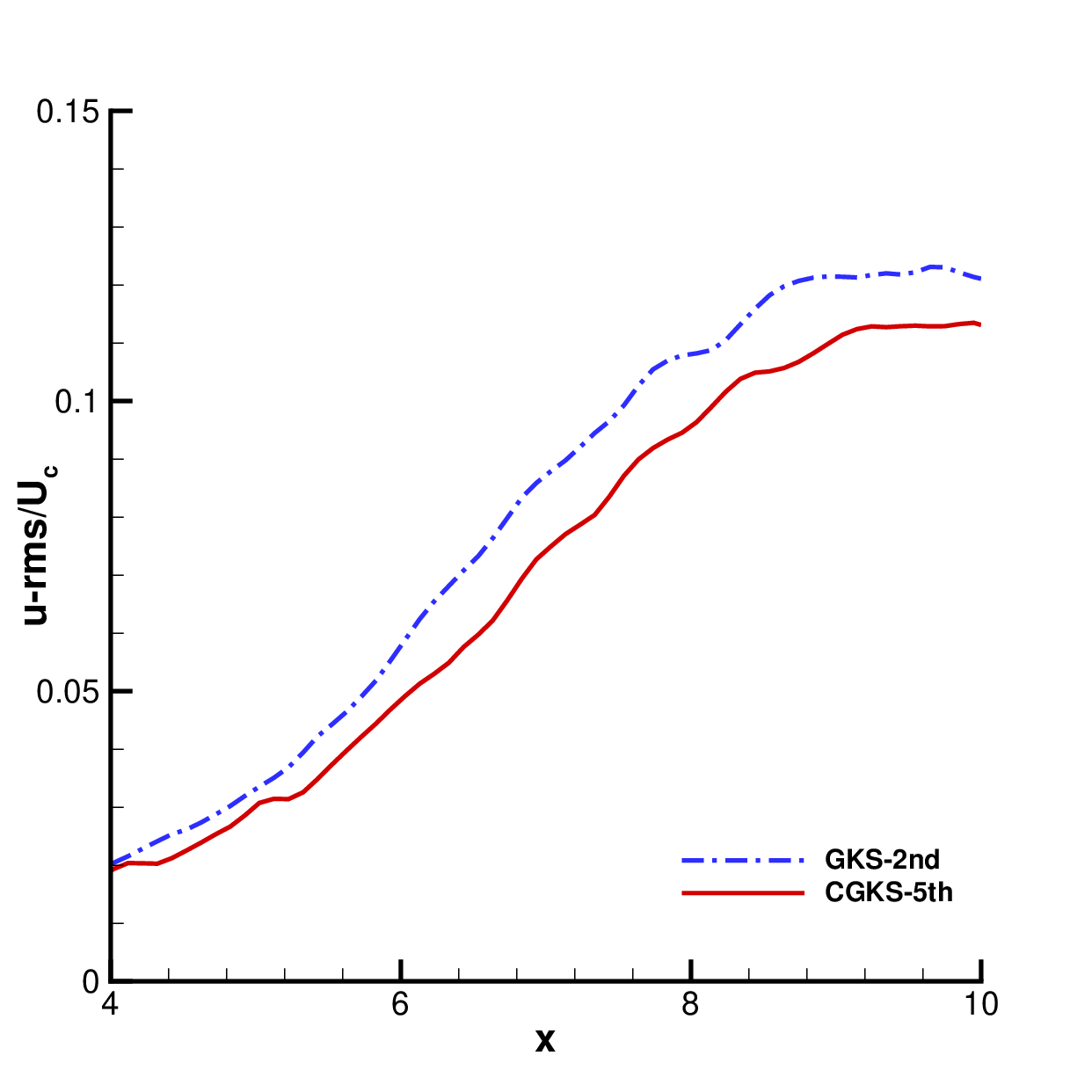}
\caption{\label{jet-sup-ave2} Supersonic Jet Simulation: the root mean square axial velocity (left) and local magnification (right) along the centerline.}
\end{figure}

\subsection{Transonic and Supersonic Jet Simulation}

Jet simulation plays a critical role in advancing the understanding of jet flows and optimizing them for a wide range of engineering applications. In transonic and supersonic jet simulations, computation is particularly challenging due to the involvement of complex physical phenomena, such as the interaction of shear layers and vortex structures, as well as the transition from laminar to turbulent flow.  
To accurately and efficiently capture these intricate dynamics, a combination of high-order numerical schemes and large-scale computations is typically required to perform long-duration, scale-resolving simulations.
As discussed in the previous subsection, achieving comparable resolution with a second-order scheme may require tens of times more cells than a fifth-order scheme, making the simulation prohibitively expensive at that scale.
In this subsection, the accuracy and resolution of CGKS-5th and GKS-2nd are compared under the constraint of equal computational resources. Specifically, the performance of the schemes is evaluated at matched computational cost by employing different mesh sizes, with mesh distributions carefully designed to ensure a fair comparison.

For this case, the isothermal jet with Mach numbers $Ma_{\infty} = 0.9$ and $Ma_{\infty} = 1.5$ are considered, with the Reynolds number set as $Re=4\times 10^5$ based on a jet diameter of $D =2r_0= 1$.
None of the simulations include any nozzle geometry.
The inflow axial velocity of the jet is initialized through a hyperbolic-tangent profile similarly to previous studies \cite{Jet-ref-1, Jet-ref-2}
\begin{equation*}
u(r)=\displaystyle\frac{U_c}{2}+\frac{U_c}{2}\tanh (\frac{r_0-r}{2\delta_\theta(0)}),
\end{equation*}
where $U_c=Ma_{\infty}\cdot c_\infty$ is the inflow centerline velocity, $c_\infty=1$, and the initial momentum thickness $\delta_\theta(0)$ is set as $0.03D$.
The inlet temperature profile is calculated using the Crocco-Busemann relation, assuming constant pressure
\begin{equation*}
\displaystyle\frac{T}{T_c}=\frac{T_\infty}{T_c}+\big[1-\frac{T_\infty}{T_c}+\frac{\gamma-1}{2}{Ma_{\infty}}^2\big(1-\sigma \frac{u(r)}{U_c}\big)\big]\frac{u(r)}{U_c},
\end{equation*}
where $T_c$ and $T_\infty$ are the jet centerline and ambient temperatures, respectively, in the present isothermal case, $T_c=T_\infty$ and $\sigma={Pr}^{1/3}$ is the recovery factor, the Prandtl number $Pr$ is set to 0.7.
Outside the inlet, the initial condition is given as follows,
\begin{align*}
(\rho,U,V,W,p)_c = (\gamma, 0, 0, 0, 1).
\end{align*}
No disturbance is added to the field.
The non-reflecting boundary condition and sponge zone are implemented in the outflow directions to filter out possible reflected waves.

In this case, a structured mesh consisting of $352\times 176\times 176$ cells is designed for CGKS-5th, as illustrated in Figure \ref{jet-mesh}.
In the flow direction, the mesh size is constant at $\Delta x=0.04D$ up to $x=9.6D$ to resolve the domain near the nozzle exit. Between $9.6D < x < 15D$, the mesh size increases exponentially from $0.04D$ to $0.16D$. Beyond this, a stretched mesh forms a sponge zone with $\Delta x_{max}=1.6D$ to damp vortices and perturbations.
Radially, the mesh size is $0.025D$ at the jet center, decreasing exponentially to $0.015D$ at $r=r_0$ to resolve small-scale structures near the jet boundary. Outside the jet, the mesh size increases back to $0.025D$ and gradually expands to $0.2D$ within $r<4.5D$. Beyond this, a sponge region is created with a maximum mesh size of $0.4D$.
To achieve computational cost equivalence with the fifth-order scheme using the above mesh, while maintaining a consistent stretching design, a structured mesh consisting of $528 \times 264 \times 264$ cells is designed for GKS-2nd, where the size of each cell is approximately $1/1.5$ times the cell size of the original mesh at the corresponding locations.
Table \ref{jet-ratio} presents the computational time required for CGKS-5th on a mesh consisting of $352 \times 176 \times 176$ cells and GKS-2nd on a mesh consisting of $528 \times 264 \times 264$ cells. It can be observed that the computational costs of the two schemes are approximately equivalent, with a cost difference within 5$\%$.
\begin{table}[!h]
\begin{center}
\def\temptablewidth{0.9\textwidth}{\rule{\temptablewidth}{1.0pt}}
\begin{tabular*}{\temptablewidth}{@{\extracolsep{\fill}}c|c|c|c|c} 
Number of GPUs & MESH  & Scheme & Time & Ratio  \\
\hline
8 & $352 \times 176 \times 176$ & CGKS-5th & 7502.3 & 1 \\
8 & $528 \times 264 \times 264$ & GKS-2nd  & 7115.9 & 0.95 \\
\end{tabular*}
{\rule{\temptablewidth}{1.0pt}}
\end{center}
\caption{\label{jet-ratio} Turbulent Flow Past a Cylinder: the evaluation of the computational efficiency of CGKS-5th and GKS-2nd for $t \in [0,10]$ with a computational cost difference within 5$\%$.}
\end{table}

For the transonic jet simulation with $Ma_{\infty} = 0.9$, Figures \ref{jet-sub-c5-Q} and \ref{jet-sub-p1-Q} provide instantaneous visualizations of the iso-surfaces of the Q-criterion, which are colored by the Mach number, alongside grayscale pressure maps illustrating pressure fluctuations relative to $P_0$. 
The results are obtained with CGKS-5th and GKS-2nd, respectively. The higher-order CGKS-5th clearly resolves more intricate turbulent structures than the second-order GKS-2nd.
Figure \ref{jet-sub-ave1} and Figure \ref{jet-sub-ave2} present the mean axial velocity and root mean square (RMS) axial velocity along the centerline ($\boldsymbol{r}=0$). 
These results are computed using CGKS-5th on a mesh with $352 \times 176 \times 176$ cells and GKS-2nd on a mesh with $528 \times 264 \times 264$ cells. 
Theses figures compare the numerical results with experimental results from Arakeri et al.\cite{Jet-ref-4} and Power et al.\cite{Jet-ref-3}, and with results obtained using the LES method combined with a sixth-order compact scheme on a multiblock structured grid containing 24.5 million cells (Fosso Pouangu{\'e} et al.) \cite{Jet-ref-2}.
Furthermore, in Figure \ref{jet-sub-ave2}, the position is normalized by the corresponding potential core length, $L_C$, defined as the distance from the jet inlet to the point where the centerline velocity drops to 95$\%$ of the inflow centerline velocity. 
The potential core length values obtained from CGKS-5th, GKS-2nd, and the reference data are presented in Table \ref{jet-Lc}. 
As shown, $L_C$ computed by CGKS-5th agrees more closely with the reference data than that from GKS-2nd.
As shown in Figure \ref{jet-sub-ave1} and Figure \ref{jet-sub-ave2}, the mean axial velocity computed using CGKS-5th demonstrates excellent agreement with both the experimental data \cite{Jet-ref-4, Jet-ref-3} and the numerical results obtained from the sixth-order LES on a refined mesh \cite{Jet-ref-2}. In contrast, GKS-2nd exhibits significant discrepancies. Similarly, the RMS axial velocity predicted by CGKS-5th also aligns well with the LES results \cite{Jet-ref-2} and the experimental data from Arakeri et al. \cite{Jet-ref-4}, whereas the second-order scheme displays clear deviations.
Notably, in the region $x \in [4,10]$, where the turbulent structures are most complex, the results computed by GKS-2nd show pronounced discrepancies in both Figure \ref{jet-sub-ave1} and Figure \ref{jet-sub-ave2}. These results demonstrate that, at matched computational cost, CGKS-5th delivers higher accuracy and resolution, yielding more reliable turbulence predictions.

\begin{table}[!h]
\begin{center}
\def\temptablewidth{0.5\textwidth}{\rule{\temptablewidth}{1.0pt}}
\begin{tabular*}{\temptablewidth}{@{\extracolsep{\fill}}c|c} 
Method & $L_C$  \\
\hline
Exp. Arakeri et al. \cite{Jet-ref-4} & 7.7D \\
Exp. Power et al. \cite{Jet-ref-3} & 7.0D \\
LES Fosso Pouangu{\'e} et al. \cite{Jet-ref-2} & 6.2D \\
GKS-2nd & 5.6D \\
CGKS-5th & 6.6D \\
\end{tabular*}
{\rule{\temptablewidth}{1.0pt}}
\end{center}
\caption{\label{jet-Lc} Transonic Jet Simulation: the comparison of $L_C$ calculated by CGKS-5th and GKS-2nd with reference data.}
\end{table}

For the supersonic jet simulation with $Ma_{\infty} = 1.5$, Figure \ref{jet-sup-c5-Q} and Figure \ref{jet-sup-p1-Q} show instantaneous visualizations of the iso-surfaces of the Q-criterion, colored by the Mach number, along with grayscale pressure visualizations representing pressure fluctuations around $P_0$, obtained using CGKS-5th and GKS-2nd, respectively. The higher-order CGKS-5th again demonstrates superior capability in resolving flow-field details.
Figure \ref{jet-sup-shadow} presents the instantaneous shadowgraph contours obtained using CGKS-5th and GKS-2nd. Both results show the propagation of sound waves within the computational domain. Notably, CGKS-5th results reveal clearer small-scale structures and sharper internal shock features within the jet, evidencing the higher effective resolution of the high-order method.
Figure \ref{jet-sup-ave1} and Figure \ref{jet-sup-ave2} compare the mean axial velocity and root-mean-square (RMS) axial velocity obtained using CGKS-5th and GKS-2nd.
The oscillations in the mean axial velocity computed by CGKS-5th over $x \in [0,5]$ align with the locations of internal shock waves identified in Figure \ref{jet-sup-shadow}, whereas the second-order GKS-2nd fails to capture this signature.
These results demonstrate that, at matched computational cost, CGKS-5th achieves substantially higher resolution for complex turbulent flows featuring shocks and small-scale vortical structures, underscoring its practical value for turbulence simulations.

\section{Conclusion}

This study presents a performance comparison between CGKS-5th and GKS-2nd on structured meshes for complex viscous flows, with a focus on computational efficiency and simulation accuracy.  Two comparative approaches are employed: (i) the computational cost  required to achieve comparable resolution, and (ii) the accuracy and resolution attainable under matched computational resources.
Across both smooth and discontinuous flows, CGKS-5th achieves the same resolution at approximately a 7-9 times lower computational cost than GKS-2nd. Under equal resource constraints, CGKS-5th further delivers substantially higher accuracy and resolution, with particularly pronounced benefits for turbulent flows involving shocks and small-scale vortical structures.
The advantages of high-order schemes for smooth compressible flows have been well-recognized in previous studies. However, their effectiveness in complex compressible flow simulations, where nonlinear reconstructions are introduced to handle discontinuities and may reduce the benefits of high-order schemes, has lacked thorough validation. For the first time, this study demonstrates, through simulations of turbulence benchmarks and jet flow problems, that high-order schemes, particularly compact schemes with high-resolution properties, offer significant advantages in both accuracy and efficiency over second-order schemes, underscoring their potential value for engineering applications.

\section*{Acknowledgements}

The current research is supported by National Key R$\&$D Program of China (Grant Nos.2022YFA1004500), National Science Foundation of China (92371107, 12172316), and Hong Kong research grant council (16301222, 16208324).

\section*{Declaration of competing interest}
The authors declare that they have no known competing financial interests or personal relationships 
that could have appeared to influence the work reported in this paper.

\section*{Data availability}
The data that support the findings of this study are available from
the corresponding author upon reasonable request.

\end{document}